# NECESSARY AND SUFFICIENT LYAPUNOV-LIKE CONDITIONS FOR ROBUST NONLINEAR STABILIZATION

**Iasson Karafyllis[*] and Zhong-Ping Jiang[**]**

**[*]Department of Environmental Engineering, Technical University of Crete, 73100, Chania, Greece**
**email: ikarafyl@enveng.tuc.gr**

**[**]Department of Electrical and Computer Engineering, Polytechnic University, Six Metrotech Center, Brooklyn, NY 11201, U.S.A.**
**email: zjiang@control.poly.edu**

## Abstract

In this work, we propose a methodology for the expression of necessary and sufficient Lyapunov-like conditions for the existence of stabilizing feedback laws. The methodology is an extension of the well-known Control Lyapunov Function (CLF) method and can be applied to very general nonlinear time-varying systems with disturbance and control inputs, including both finite- and infinite-dimensional systems. The generality of the proposed methodology is also reflected upon by the fact that partial stability with respect to output variables is addressed. In addition, it is shown that the generalized CLF method can lead to a novel tool for the explicit design of robust nonlinear controllers for a class of time-delay nonlinear systems with a triangular structure.

**Keywords:** Control Lyapunov Function, stabilization, time-varying systems, nonlinear control.

## 1. Introduction

Feedback stabilization of nonlinear systems is a fundamentally important problem in control theory and practice. The purpose of this paper is to look at this problem from a control Lyapunov function point of view, but for a wide class of nonlinear time-varying systems. We aim to develop a methodology that not only results in necessary and sufficient conditions for robust feedback stabilization, but provides novel tools for the design of robust nonlinear controllers. To add to the generality of this framework, we will address partial stability with respect to output variables, instead of state variables. We first consider finite-dimensional nonlinear systems, and then show that the same methodology can be adapted to infinite-dimensional systems described by retarded functional differential equations.

Specifically, we begin with finite-dimensional control systems in the general form:

$$\begin{aligned} &\dot{x}(t) = f(t, d(t), x(t), u(t)) \\ &Y(t) = H(t, x(t)) \\ &x(t) \in \Re^n \,, d(t) \in D \,, t \geq 0 \,, Y(t) \in \Re^k \,, u(t) \in U \subseteq \Re^m \end{aligned} \tag{1.1}$$

where the vector fields $f : \Re^+ \times D \times \Re^n \times U \to \Re^n$, $H : \Re^+ \times \Re^n \to \Re^k$ are continuous with $f(t,d,0,0) = 0$, $H(t,0) = 0$ for all $(t,d) \in \Re^+ \times D$. And we ask the following question of feedback stabilizability: Under what conditions there exists a continuous feedback of the form:

$$u = k(t, x) \tag{1.2}$$

such that the closed-loop system (1.1) with (1.2) is (uniformly) Robustly Globally Asymptotically Output Stable? See Section 2.1 for a precise definition.

The above-mentioned problem has been studied by several authors in past literature for a subclass of nonlinear control systems (1.1). For instance, in his pioneering work [1] Artstein studied the above existence problem for affine autonomous control systems without disturbances, $U \subseteq \Re^m$ being a closed convex set and output $Y$ being identically the state of the system, i.e., $H(t,x) \equiv x$ (see also [27]). He showed in [1] that the existence of a time-independent Control Lyapunov Function (CLF) satisfying the "small-control" property is a necessary and sufficient condition for the existence of a continuous stabilizing feedback. Sontag [24] extended the results by presenting an explicit formula of the feedback stabilizer for affine autonomous control systems without disturbances, $U = \Re^m$ and output $Y$ being identically the state of the system. Sontag's formula was exploited recently in [9] for the uniform stabilization of time-varying systems. Freeman and Kokotovic in [4] extended the idea of the CLF in order to study affine control systems with disturbances, $U \subseteq \Re^m$ being a closed convex set and output $Y$ being identically the state of the system, i.e., $H(t,x) \equiv x$ : they introduced the concept of the Robust Control Lyapunov Function (RCLF). In [10] the authors showed that the "small-control" property is not needed for non-uniform in time robust global stabilization of the state ($H(t,x) \equiv x$ ) of control systems affine in the control with $U = \Re^m$ . The result was extended in [13] for the general case of output stability. In all the above approaches the stabilizing feedback is constructed using a partition of unity methodology or Michael's Theorem (when simple continuity of the feedback suffices). Control Lyapunov Functions have also been used for the design of discontinuous feedback laws (see for instance [3]), the design of static output feedback stabilizers (see [13,28]), as well as for the design of adaptive nonlinear controllers (see [17,23]).

However, so far the method of Lyapunov design of stabilizing feedback laws is more frequently applied to finite-dimensional systems of the form (1.1). In order to be able to extend the applicability of the method to infinite-dimensional systems of the form $\dot{x} = f(t,d,x,u)$ where the state $x$ belongs to an infinite-dimensional normed linear space $\mathcal{X}$, one has to deal with Control Lyapunov Functionals $V : \Re^+ \times \mathcal{X} \to \Re^+$ , which present one (or many) of the following complications:

(*i*) In contrast to CLF in the finite-dimensional case, usually Control Lyapunov Functionals are simply locally Lipschitz mappings of the state (and, not necessarily, continuously differentiable);

(*ii*) Even if the mapping $f$ is affine in $u$, the (appropriate) derivative of the Control Lyapunov Functional $\dot{V}(t,d,x,u)$ is not necessarily affine in $u$ ;

(*iii*) The existing feedback construction methodology based on partition of unity arguments (see. e.g., [1,27]) does not work because the state space $\mathcal{X}$ is infinite-dimensional;

(*iv*) The feedback construction methodology based on Michael's Theorem (see, e.g., [4]) does not work either because simple continuity of the feedback does not suffice or because the hypotheses of Michael's Theorem cannot be verified.

Particularly, all of the above complications are encountered when control systems described by Retarded Functional Differential Equations (RFDEs) are studied, i.e., systems of the form

$$\begin{aligned} &\dot{x}(t) = f(t, d(t), T_r(t)x, u(t)) \,,\; t \ge t_0 \\ &Y(t) = H(t, T_r(t)x) \\ &x(t) \in \Re^n \;, Y(t) \in \mathcal{Y},\; d(t) \in D \,, u(t) \in U \end{aligned} \tag{1.3}$$

where $r > 0$ is a constant, $f : \Re^+ \times D \times C^0([-r,0];\Re^n) \times U \to \Re^n$, $H : \Re^+ \times C^0([-r,0];\Re^n) \to \mathcal{Y}$ satisfy $f(t,d,0,0) = 0$, $H(t,0) = 0$ for all $(t,d) \in \Re^+ \times D$, $D \subseteq \Re^l$ is a non-empty compact set, $U \subseteq \Re^m$ is a closed convex set with $0 \in U$ , $\mathcal{Y}$ is a normed linear space and $T_r(t)x = x(t+\theta)$ ; $\theta \in [-r,0]$. It should be emphasized that by allowing the output to take values in an abstract normed linear spaces we are in a position to consider:

- outputs with no delays, e.g. $Y(t) = h(t, x(t))$ with $\mathcal{Y} = \Re^k$ ,
- outputs with discrete or distributed delay, e.g. $Y(t) = h(x(t), x(t-r))$ or $Y(t) = \int_{t-r}^{t} h(t,\theta,x(\theta))d\theta$ with $\mathcal{Y} = \Re^k$ ,
- functional outputs with memory, e.g. $Y(t) = h(t,\theta,x(t+\theta))$ ; $\theta \in [-r,0]$ or the identity output $Y(t) = T_r(t)x = x(t+\theta)$ ; $\theta \in [-r,0]$ with $\mathcal{Y} = C^0([-r,0];\Re^k)$ .

As the second contribution of the present work, we show how all complications mentioned above for infinite-dimensional systems can be solved, and consequently we obtain Lyapunov-like necessary and sufficient conditions for systems of the form (1.3). Since the methodology that we describe in the present work allows the construction of locally Lipschitz stabilizing feedback laws, it is expected that it can be used for general infinite-dimensional control systems. More importantly, we will show that our generalized CLF methodology is more than of existence-type result, but can yield constructive design tools for an enlarged class of nonlinear control systems. To this end, we will study in details a class of triangular time-delay nonlinear systems described by RFDEs, i.e.

$$\begin{aligned}
&\dot{x}_i(t) = f_i(t,d(t),T_r(t)x_1,...,T_r(t)x_i) + g_i(t,d(t),T_r(t)x_1,...,T_r(t)x_i)x_{i+1}(t) \quad , \quad i=1,...,n-1\\
&\dot{x}_n(t) = f_n(t,d(t),T_r(t)x) + g_n(t,d(t),T_r(t)x)u(t)\\
&x(t)=(x_1(t),...,x_n(t))\in \Re^n \; , d(t)\in D \, , u(t)\in \Re \, , t\ge 0
\end{aligned} \tag{1.4}$$

Autonomous and disturbance-free systems of the form (1.4) have been studied in [7,8,22,29]. In the present work it is shown that the construction of a stabilizing feedback law for (1.4) proceeds in parallel with the construction of a State Robust Control Lyapunov Functional. Moreover, sufficient conditions for the existence and design of a stabilizing feedback law $u(t)=k(x(t))$, which is independent of the delay are given.

The rest of the paper is organized as follows. First, our results are developed for the finite-dimensional case (1.1) (Section 2), where necessary and sufficient conditions for the existence of stabilizing feedback are formulated. Section 3 is devoted to the development of sufficient conditions, which guarantee that a given function is an Output Robust Control Lyapunov function for (1.1). Examples are presented for systems which are polynomial in the control; this case was recently studied in [21]. In Section 4 we show how the same methodology can be applied to the infinite-dimensional case (1.3). We exploit the converse Lyapunov theorems in [14,16] to obtain necessary and sufficient conditions for the uncertain case (1.3); the classical Lyapunov-Krasovskii characterizations (see [5]) cannot be used since they are applicable to the disturbance-free case. Lyapunov-based feedback design for systems described by RFDEs was used recently in [6,7,19,20,22,29]. Section 5 is devoted to the case of triangular control systems (1.4). Finally, our concluding remarks are given in Section 6.

**Notations** Throughout this paper we adopt the following notations:

* Let $A\subseteq \Re^n$ be a set. By $C^0(A;\Omega)$, we denote the class of continuous functions on $A$, which take values in $\Omega$. By $C^k(A;\Omega)$, where $k\ge 1$ is an integer, we denote the class of differentiable functions on $A$ with continuous derivatives up to order $k$, which take values in $\Omega$. By $C^\infty(A;\Omega)$, we denote the class of differentiable functions on $A$ having continuous derivatives of all orders, which take values in $\Omega$, i.e., $C^\infty(A;\Omega)=\bigcap_{k\ge 1} C^k(A;\Omega)$.
* By $\| \ \|_Y$, we denote the norm of the normed linear space $Y$.
* A continuous mapping $A\times B \ni (z,x)\to k(z,x)\in\Re^m$, where $B\subseteq X$, $A\subseteq Y$ and $X,Y$ are normed linear spaces, is called completely locally Lipschitz with respect to $x\in B$ if for every closed and bounded set $S\subseteq A\times B$ it holds that $\sup\left\{ \frac{|k(z,x)-k(z,y)|}{\|x-y\|_X} : (z,x)\in S\, ,(z,y)\in S\, , x\ne y \right\} < +\infty$. If the normed linear spaces $X,Y$ are finite-dimensional spaces then we simply say that the continuous mapping $A\times B\ni (z,x)\to k(z,x)\in \Re^m$ is locally Lipschitz with respect to $x\in B$ if for every compact set $S\subseteq A\times B$ it holds that $\sup\left\{ \frac{|k(z,x)-k(z,y)|}{|x-y|} : (z,x)\in S\, ,(z,y)\in S\, , x\ne y \right\} < +\infty$.
* For a vector $x\in\Re^n$ we denote by $|x|$ its usual Euclidean norm and by $x'$ its transpose. For $x\in C^0([-r,0];\Re^n)$ we define $\|x\|_r := \max_{\theta\in[-r,0]} |x(\theta)|$.
* $\Re^+$ denotes the set of non-negative real numbers.
* $\mathcal{E}$ denotes the class of non-negative $C^0$ functions $\mu:\Re^+\to\Re^+$, for which it holds: $\int_0^{+\infty}\mu(t)dt<+\infty$ and $\lim_{t\to+\infty}\mu(t)=0$.
* We denote by $K^+$ the class of positive $C^0$ functions defined on $\Re^+$. We say that a function $\rho:\Re^+\to\Re^+$ is

positive definite if $\rho(0)=0$ and $\rho(s)>0$ for all $s>0$. By $K$ we denote the set of positive definite, increasing and continuous functions. We say that a positive definite, increasing and continuous function $\rho:\Re^+\to\Re^+$ is of class $K_\infty$ if $\lim_{s\to+\infty}\rho(s)=+\infty$. By $KL$ we denote the set of all continuous functions $\sigma=\sigma(s,t):\Re^+\times\Re^+\to\Re^+$ with the properties: (i) for each $t\geq 0$ the mapping $\sigma(\cdot,t)$ is of class $K$ ; (ii) for each $s\geq 0$, the mapping $\sigma(s,\cdot)$ is non-increasing with $\lim_{t\to+\infty}\sigma(s,t)=0$.

* Let $D\subseteq\Re^l$ be a non-empty set. By $M_D$ we denote the class of all Lebesgue measurable and locally essentially bounded mappings $d:\Re^+\to D$.
* Let $x:[a-r,b)\to\Re^n$ be a continuous mapping with $b>a>-\infty$ and $r>0$. By $T_r(t)x$ we denote the "$r$-history" of $x$ at time $t\in[a,b)$, i.e., $T_r(t)x:=x(t+\theta)$ ; $\theta\in[-r,0]$. Notice that $T_r(t)x\in C^0([-r,0];\Re^n)$.
* A function $\Psi:A\times U\to\Re\cup\{+\infty\}$, where $A\subseteq\Re^n$, $U\subseteq\Re^m$ is a convex set, is called quasi-convex with respect to $u\in U$, if for every $u,v\in U$, $x\in A$ and $\lambda\in[0,1]$ it holds that $\Psi(x,\lambda u+(1-\lambda)v)\leq\max\left\{\Psi(x,u),\Psi(x,v)\right\}$.
* Let $A\subseteq\Re^n$ be a non-empty set. By $\overline{A}$ we denote the closure of $A\subseteq\Re^n$ and by $\overline{coA}$, we denote the closure of the convex hull of $A\subseteq\Re^n$.
* Let $U\subseteq\Re^n$ be a non-empty closed set. By $\mathrm{Pr}_U(x)$, we denote the projection of $x\in\Re^n$ on $U\subseteq\Re^n$. Notice that if $U\subseteq\Re^n$ is convex then $\left|\mathrm{Pr}_U(x)-\mathrm{Pr}_U(y)\right|\leq|x-y|$, for all $x,y\in\Re^n$.

## 2. Finite-Dimensional Control Systems

In this section, we consider control systems of the form (1.1) under the following hypotheses:

**(H1)** the vector fields $f:\Re^+\times D\times\Re^n\times U\to\Re^n$, $H:\Re^+\times\Re^n\to\Re^k$ are continuous and for every bounded interval $I\subset\Re^+$ and every compact set $S\subset\Re^n\times U$ there exists $L\geq 0$ such that $\left|f(t,d,x,u)-f(t,d,y,v)\right|\leq L|x-y|+L|u-v|$ for all $(t,d)\in I\times D$, $(x,u)\in S$, $(y,v)\in S$ (i.e., the mapping $\Re^+\times D\times\Re^n\times U\ni(t,d,x,u)\to f(t,d,x,u)\in\Re^n$ is locally Lipschitz with respect to $(x,u)$),

**(H2)** the set $D\subset\Re^l$ is compact and $U\subseteq\Re^m$ is a closed convex set,

**(H3)** $f(t,d,0,0)=0$, $H(t,0)=0$ for all $(t,d)\in\Re^+\times D$.

In order to present the main results on finite-dimensional systems of the form (1.1) we need to present in detail the basic steps of the method. The methodology consists of the following steps:

**2.I.** Notions of Output Stability
**2.II.** Lyapunov-like criteria for Output stability
**2.III.** Definition of the Output Robust Control Lyapunov Function
**2.IV.** Converse Lyapunov theorems for output stability

### 2.I. Notions of Output Stability

We first analyze the output stability notions used in the present work. Consider the system

$$\begin{aligned}&\dot{x}(t)=f(t,d(t),x(t))\\&Y(t)=H(t,x(t))\\&x(t)\in\Re^n\,,d(t)\in D\,,Y(t)\in\Re^k\end{aligned}\tag{2.1}$$

where the vector fields $f:\Re^+\times D\times\Re^n\to\Re^n$, $H:\Re^+\times\Re^n\to\Re^k$ are continuous and $D\subset\Re^l$ is compact. We assume that for every $(t_0, x_0, d) \in \Re^+ \times \Re^n \times M_D$ there exist $h \in (0,+\infty]$ and a unique absolutely continuous mapping $x : [t_0, t_0 + h) \to \Re^n$ with $x(t_0) = x_0$ and $\dot{x}(t) = f(t, d(t), x(t))$ a.e. for $t \in [t_0, t_0 + h)$. Moreover, we assume that $f(t, d, 0) = 0$, $H(t,0) = 0$ for all $(t, d) \in \Re^+ \times D$. The solution $x : [t_0, t_0 + h) \to \Re^n$ of (2.1) at time $t \geq t_0$ with initial condition $x(t_0) = x_0$ corresponding to input $d \in M_D$ will be denoted by $x(t, t_0, x_0; d)$.

**Definition 2.1:** *We say that (2.1) is **Robustly Forward Complete (RFC)** if for every $T \geq 0$, $r \geq 0$ it holds that:*

$$\sup\left\{ \left| x(t_0 + h, t_0, x_0; d) \right| ;\ \left| x_0 \right| \leq r \,,\ t_0 \in [0,T] \,,\ h \in [0,T] \,,\ d(\cdot) \in M_D \right\} < +\infty \tag{2.2}$$

Clearly, the notion of robust forward completeness implies the standard notion of forward completeness, which simply requires that for every initial condition the solution of the system exists for all times greater than the initial time, or equivalently, the solutions of the system do not present finite escape time. Conversely, an extension of Proposition 5.1 in [18] to the time-varying case shows that every forward complete system (2.1) whose dynamics are locally Lipschitz with respect to $(t, x)$, uniformly in $d \in D$, is RFC. All output stability notions used in the present work will assume RFC.

We continue with the notion of (non-uniform in time) Robust Global Asymptotic Output Stability (RGAOS) as a generalization of the notion of Robust Output Stability (see [11,12]). Let us denote by $Y(t) = H(t, x(t, t_0, x_0; d))$ the value of the output for the unique solution of (2.1) at time $t$ that corresponds to input $d \in M_D$ with initial condition $x(t_0) = x_0$.

**Definition 2.2:** *Consider system (2.1) and suppose that (2.1) is RFC. We say that system (2.1) is **(non-uniformly in time) Robustly Globally Asymptotically Output Stable (RGAOS)** if it satisfies the following properties:*

***P1(Output Stability)*** *For every $\varepsilon > 0$, $T \geq 0$, it holds that*

$$\sup\left\{ \left| Y(t) \right| ;\ t \geq t_0 \,,\ \left| x_0 \right| \leq \varepsilon \,,\ t_0 \in [0,T] \,,\ d(\cdot) \in M_D \right\} < +\infty$$

*and there exists a $\delta := \delta(\varepsilon, T) > 0$ such that:*

$$\left| x_0 \right| \leq \delta \,,\ t_0 \in [0,T] \Rightarrow \left| Y(t) \right| \leq \varepsilon \,,\ \forall t \geq t_0 \,,\ \forall d(\cdot) \in M_D$$

***P2(Uniform Output Attractivity on compact sets of initial data)*** *For every $\varepsilon > 0$, $T \geq 0$ and $R \geq 0$, there exists a $\tau := \tau(\varepsilon, T, R) \geq 0$, such that:*

$$\left| x_0 \right| \leq R \,,\ t_0 \in [0,T] \Rightarrow \left| Y(t) \right| \leq \varepsilon \,,\ \forall t \geq t_0 + \tau \,,\ \forall d(\cdot) \in M_D$$

The notion of Uniform Robust Global Asymptotic Output Stability was given in [25,26] and is a special case of (non-uniform in time) RGAOS.

**Definition 2.3:** *Consider system (2.1) and suppose that (2.1) is RFC. We say that system (2.1) is **Uniformly Robustly Globally Asymptotically Output Stable (URGAOS)** if it satisfies the following properties:*

***P1(Uniform Output Stability)*** *For every $\varepsilon > 0$, it holds that*

$$\sup\left\{ \left| Y(t) \right| ;\ t \geq t_0 \,,\ \left| x_0 \right| \leq \varepsilon \,,\ t_0 \geq 0 \,,\ d(\cdot) \in M_D \right\} < +\infty$$

*and there exists a $\delta := \delta(\varepsilon) > 0$ such that:*

$$\left| x_0 \right| \leq \delta \,,\ t_0 \geq 0 \Rightarrow \left| Y(t) \right| \leq \varepsilon \,,\ \forall t \geq t_0 \,,\ \forall d(\cdot) \in M_D$$

***P2(Uniform Output Attractivity on compact sets of initial states)*** *For every* $\varepsilon > 0$ *and* $R \ge 0$, *there exists a* $\tau := \tau(\varepsilon, R) \ge 0$, *such that:*

$$|x_0| \le R,\, t_0 \ge 0 \Rightarrow |Y(t)| \le \varepsilon,\, \forall t \ge t_0 + \tau,\ \forall d(\cdot) \in M_D$$

Obviously, for the case $H(t,x) = x$ the notions of RGAOS and URGAOS coincide with the notions of ***non-uniform in time Robust Global Asymptotic Stability (RGAS)*** as given in [10] and ***Uniform Robust Global Asymptotic Stability (URGAS)*** as given in [18], respectively. Also note that if there exists $a \in K_\infty$ with $|x| \le a(|H(t,x)|)$ for all $(t,x) \in \Re^+ \times \Re^n$, then (U)RGAOS implies (U)RGAS .

## 2.II. Lyapunov-like Criteria for Output Stability

For a locally bounded function $V : \Re^+ \times \Re^n \to \Re$, we define

$$V^0(t,x;v) := \limsup_{h \to 0^+,\, w \to v} \frac{V(t+h, x+hw) - V(t,x)}{h} \tag{2.3}$$

The reader should notice that the function $(t,x,v) \to V^0(t,x;v)$ may take values in the extended real number set $\Re^* = [-\infty, +\infty]$. However, for locally Lipschitz functions $V : \Re^+ \times \Re^n \to \Re$, the function $(t,x,v) \to V^0(t,x;v)$ is locally bounded. It should be clear that for locally Lipschitz functions $V : \Re^+ \times \Re^n \to \Re$ it holds that:

$$V^0(t,x;v) = \limsup_{h \to 0^+} \frac{V(t+h, x+hv) - V(t,x)}{h} \tag{2.4}$$

The main reason for introducing the above Dini derivative is the following lemma.

**Lemma 2.4:** *Let* $V : \Re^+ \times \Re^n \to \Re$ *be a locally bounded function and let* $x : [t_0, t_{\max}) \to \Re^n$ *be a solution of (2.1) with initial condition* $x(t_0) = x_0 \in \Re^n$ *corresponding to certain* $d \in M_D$, *where* $t_{\max} \in (t_0, +\infty]$ *is the maximal existence time of the solution. Then it holds that*

$$\limsup_{h \to 0^+} h^{-1}\big(V(t+h, x(t+h)) - V(t, x(t))\big) \le V^0(t, x; D^+x(t)),\ \text{a.e. on } [t_0, t_{\max}) \tag{2.5}$$

*where* $D^+x(t) = \lim_{h \to 0^+} h^{-1}\big(x(t+h) - x(t)\big)$.

**Proof** Inequality (2.5) follows directly from definition (2.3) and definition $w_h = \dfrac{x(t+h) - x(t)}{h}$, where $t \in [t_0, t_{\max}) \setminus N$ and $N$ is a measure zero set where $D^+x(t) = \lim_{h \to 0^+} h^{-1}\big(x(t+h) - x(t)\big)$ is not defined. Notice that since $D^+x(t) = \lim_{h \to 0^+} h^{-1}\big(x(t+h) - x(t)\big)$ we obtain that $w_h \to D^+x(t)$ as $h \to 0^+$. The proof is complete. $\lhd$

Having introduced an appropriate derivative for Lyapunov functions, we are now in a position to give Lyapunov-like criteria for RGAOS and URGAOS. The proof of the following proposition can be found in the Appendix.

**Proposition 2.5:** *Consider system (2.1) and the following statements:*

**(Q1)** *There exist a locally Lipschitz function* $V : \Re^+ \times \Re^n \to \Re^+$, *functions* $a_1, a_2 \in K_\infty$, $\beta, \mu \in K^+$, *a function* $q \in \mathfrak{E}$ *and a* $C^0$ *positive definite function* $\rho : \Re^+ \to \Re^+$ *such that*

$$a_1\big(|(\mu(t)x, H(t,x))|\big) \le V(t,x) \le a_2\big(\beta(t)|x|\big),\ \forall (t,x) \in \Re^+ \times \Re^n \tag{2.6}$$

*and such that the following inequality holds for all* $(t,x,d) \in \Re^+ \times \Re^n \times D$ *:*

$$V^0(t,x;\, f(t,d,x)) \le -\rho\big(V(t,x)\big) + q(t) \tag{2.7}$$

**(Q2)** *Hypothesis (Q1) holds with* $\beta(t) \equiv 1$ *,* $q(t) \equiv 0$ *.*

**(Q3)** *The mapping* $\Re^+ \times D \times \Re^n \ni (t,d,x) \to f(t,d,x) \in \Re^n$ *is locally Lipschitz with respect to* $x \in \Re^n$ *.*

*If hypotheses (Q1), (Q3) hold then system (2.1) is RGAOS. If hypothesis (Q2) holds, then system (2.1) is URGAOS.*

## 2.III. Definition of the Output Robust Control Lyapunov Function

We next give the definition of the Output Robust Control Lyapunov Function for system (1.1). The definition is in the same spirit with the definition of the notion of Robust Control Lyapunov Function given in [4] for continuous-time finite-dimensional control systems. The small-control property in the following definition constitutes a time-varying version of the small-control property for the autonomous case [1,4,24].

**Definition 2.6:** *We say that (1.1) admits an* ***Output Robust Control Lyapunov Function (ORCLF)*** *if there exists a locally Lipschitz function* $V : \Re^+ \times \Re^n \to \Re^+$ *(called the Output Robust Control Lyapunov Function), which satisfies the following properties:*

**(i)** *There exist* $a_1, a_2 \in K_\infty$ *,* $\beta, \mu \in K^+$ *such that (2.6) holds.*

**(ii)** *There exists a function* $\Psi : \Re^+ \times \Re^n \times U \to \Re \cup \{+\infty\}$ *with* $\Psi(t,0,0) = 0$ *for all* $t \ge 0$ *such that for each* $u \in U$ *the mapping* $(t,x) \to \Psi(t,x,u)$ *is upper semi-continuous, a function* $q \in \mathcal{E}$ *and a* $C^0$ *positive definite function* $\rho : \Re^+ \to \Re^+$ *such that the following inequality holds:*

$$\inf_{u \in U} \Psi(t,x,u) \le q(t) \,,\ \forall (t,x) \in \Re^+ \times \Re^n \tag{2.8}$$

*Moreover, for every finite set* $\{u_1, u_2, ..., u_p\} \subset U$ *and for every* $\lambda_i \in [0,1]$ *(* $i = 1,..., p$ *) with* $\sum_{i=1}^{p} \lambda_i = 1$*, it holds that:*

$$\sup_{d \in D} V^0\left(t,x;\, f\left(t,d,x,\sum_{i=1}^{p} \lambda_i u_i\right)\right) \le -\rho\big(V(t,x)\big) + \max\big\{\Psi(t,x,u_i)\,,\, i = 1,..., p\,\big\},\ \forall (t,x) \in \Re^+ \times \Re^n \tag{2.9}$$

*If in addition to the above there exist* $a \in K_\infty$ *,* $\gamma \in K^+$ *such that for every* $(t,x) \in \Re^+ \times \Re^n$ *there exists* $u \in U$ *with* $|u| \le a\big(\gamma(t)|x|\big)$ *such that*

$$\Psi(t,x,u) \le q(t) \tag{2.10}$$

*then we say that* $V : \Re^+ \times \Re^n \to \Re^+$ *satisfies the "small-control" property.*

*For the case* $H(t,x) \equiv x$ *we simply call* $V : \Re^+ \times \Re^n \to \Re^+$ *a State Robust Control Lyapunov Function (SRCLF).*

**Remark 2.7:** It is important to emphasize that the Dini derivative used for the Lyapunov-like criteria for output stability (Proposition 2.5) is the same derivative used in inequality (2.9) for the definition of the ORCLF.

## 2.IV. Converse Lyapunov theorems for output stability

In this section we are going to exploit the converse Lyapunov theorem for RGAOS presented in [12].

## 2.V. Main Results

We are now ready to state and prove our main results for the finite-dimensional case (1.1).

**Theorem 2.8:** *Consider system (1.1) under hypotheses (H1-3). The following statements are equivalent:*

**(a)** *There exists a $C^\infty$ function $k:\Re^+\times\Re^n\to U$ with $k(t,0)=0$ for all $t\geq 0$, in such a way that the closed-loop system (1.1) with $u=k(t,x)$ is RGAOS.*

**(b)** *There exists a $C^0$ function $k:\Re^+\times\Re^n\to U$ with $f(t,d,x,k(t,x))$ being locally Lipschitz with respect to $x$ and $f(t,d,0,k(t,0))=0$ for all $(t,d)\in\Re^+\times D$, such that the closed-loop system (1.1) with $u=k(t,x)$ is RGAOS.*

**(c)** *System (1.1) admits an ORCLF, which satisfies the small control property with $q(t)\equiv 0$.*

**(d)** *System (1.1) admits an ORCLF.*

**Theorem 2.9:** *Consider system (1.1) under hypotheses (H1-3). If system (1.1) admits an ORCLF, which satisfies the small-control property and inequalities (2.6), (2.10) with $\beta(t)\equiv 1$, $q(t)\equiv 0$, then there exists a continuous mapping $k:\Re^+\times\Re^n\to U$, with $k(t,0)=0$ for all $t\geq 0$, which is $C^\infty$ on the set $\Re^+\times(\Re^n\setminus\{0\})$, such that*

*i)* *for all $(t_0,x_0,d)\in\Re^+\times\Re^n\times M_D$ the solution $x(t)$ of the closed-loop system (1.1) with $u=k(t,x)$, i.e. the solution of*

$$\dot{x}(t)=f(t,d(t),x(t),k(t,x(t))) \tag{2.11}$$

*with initial condition $x(t_0)=x_0\in\Re^n$, corresponding to input $d\in M_D$ is unique,*

*ii)* *system (2.11) is URGAOS.*

*Moreover, if the ORCLF $V$ and the function $\Psi$ involved in property (ii) of Definition 2.6 are time independent then the continuous mapping $k$ is time invariant. Finally, if in addition there exist functions $\eta\in K^+$, $\varphi\in C^v(A;U)$ where $v\in\{1,2,...\}$, $A=\bigcup_{t\geq 0}\{t\}\times\{x\in\Re^n:|x|<4\eta(t)\}$ with $\varphi(t,0)=0$ for all $t\geq 0$, such that*

$$\Psi(t,x,\varphi(t,x))\leq 0\text{, for all }(t,x)\in\Re^+\times\Re^n\text{ with }|x|\leq 2\eta(t) \tag{2.12}$$

*then the continuous mapping $k$ is of class $C^v(\Re^+\times\Re^n;U)$.*

**Proof of Theorem 2.8:** The implications (a)$\Rightarrow$(b) and (c)$\Rightarrow$(d) are obvious and we prove implications (d)$\Rightarrow$(a), (b)$\Rightarrow$(d) and (a)$\Rightarrow$(c).

**(d)$\Rightarrow$(a)** Suppose that (1.1) admits an ORCLF. Without loss of generality, we may assume that the function $q\in\mathcal{E}$ involved in (2.8) is positive for all $t\geq 0$.

Furthermore define:

$$\Xi(t,x,u):=\Psi(t,x,u)-8q(t)\text{, }(t,x,u)\in\Re^+\times\Re^n\times U \tag{2.13a}$$

$$\Xi(t,x,u) := \Xi(0,x,u)\,,\ (t,x,u) \in (-1,0)\times\Re^n \times U \tag{2.13b}$$

The definition of $\Xi$, given by (2.13a,b), guarantees that the function $\Xi : (-1,+\infty)\times\Re^n \times U \to \Re \cup \{+\infty\}$ with $\Xi(t,0,0) = -8q(\max\{0,t\})$ for all $t > -1$ is such that, for each $u \in U$ the mapping $(t,x) \to \Xi(t,x,u)$ is upper semi-continuous. By virtue of (2.8) and upper semi-continuity of $\Xi$, it follows that for each $(t,x) \in (-1,+\infty)\times\Re^n$ there exist $u = u(t,x) \in U$ and $\delta = \delta(t,x) \in (0,t+1)$, such that

$$\Xi(\tau,y,u(t,x)) \le 0\,,\ \forall(\tau,y) \in \left\{(\tau,y) \in (-1,+\infty)\times\Re^n : |\tau - t| + |y - x| < \delta\right\} \tag{2.14}$$

Using (2.14) and standard partition of unity arguments, we can determine sequences $\{(t_i,x_i) \in (-1,+\infty)\times\Re^n\}_{i=1}^{\infty}$, $\{u_i \in U\}_{i=1}^{\infty}$, $\{\delta_i\}_{i=1}^{\infty}$ with $\delta_i = \delta(t_i,x_i) \in (0,t_i+1)$ associated with a sequence of open sets $\{\Omega_i\}_{i=1}^{\infty}$ with

$$\Omega_i \subseteq \left\{(\tau,y) \in (-1,+\infty)\times\Re^n : |\tau - t_i| + |y - x_i| < \delta_i\right\} \tag{2.15}$$

forming a locally finite open covering of $(-1,+\infty)\times\Re^n$ and in such a way that:

$$\Xi(\tau,y,u_i) \le 0\,,\ \forall(\tau,y) \in \Omega_i \tag{2.16}$$

Also, a family of smooth functions $\{\theta_i\}_{i=1}^{\infty}$ with $\theta_i(t,x) \ge 0$ for all $(t,x) \in (-1,+\infty)\times\Re^n$ can be determined with

$$supp\ \theta_i \subseteq \Omega_i \tag{2.17}$$

$$\sum_{i=1}^{\infty}\theta_i(t,x) = 1\,,\ \forall(t,x) \in (-1,+\infty)\times\Re^n \tag{2.18}$$

The facts that $\Xi(t,0,0) = -8q(t) < 0$ for all $t \ge 0$ and that the mapping $(t,x) \to \Xi(t,x,0)$ is upper semi-continuous, imply that for every $t \ge 0$ there exists $\delta(t) > 0$ such that $\Xi(\tau,y,0) \le 0$ for all $(\tau,y) \in \Re^+ \times \Re^n$ with $|\tau - t| + |y| \le \delta(t)$. Utilizing compactness of $[0,T]$ for every $T \ge 0$, we conclude that for every $T \ge 0$ there exists $\tilde{\delta}(T) > 0$ such that

$$(\tau,y) \in [0,T]\times\Re^n \text{ and } |y| \le \tilde{\delta}(T)\ \Rightarrow \Xi(\tau,y,0) \le 0 \tag{2.19}$$

Define the following function:

$$\tilde{\eta}(t) := \frac{1}{2}\left([t]+1-t\right)\tilde{\delta}\left([t]+1\right) + \frac{1}{2}\left(t-[t]\right)\tilde{\delta}\left([t]+2\right),\ t \ge 0 \tag{2.20}$$

where $[t]$ denotes the integer part of $t \ge 0$. Notice that by virtue of definition (2.20) it follows that $\tilde{\eta}(k) = \frac{1}{2}\tilde{\delta}(k+1)$, $\lim_{t \to (k+1)^-}\tilde{\eta}(t) = \frac{1}{2}\tilde{\delta}(k+2)$ for all $k \in Z^+$, which implies that $\tilde{\eta}$ is continuous. Moreover, definition (2.20) gives $0 < \tilde{\eta}(t) \le \frac{1}{2}\max\left\{\tilde{\delta}\left([t]+1\right); \tilde{\delta}\left([t]+2\right)\right\}$ for all $t \ge 0$, which in conjunction with (2.19) and the inequality $t \le [t]+1 \le [t]+2$ implies:

$$|x| \le 2\tilde{\eta}(t)\ \Rightarrow \Xi(t,x,0) \le 0 \tag{2.21}$$

Let $\bar{\eta} : \Re^+ \to (0,+\infty)$ be the positive, continuous and non-decreasing function defined by $\bar{\eta}(t) := \min_{0 \le \tau \le t}\tilde{\eta}(\tau)$. Let $\varphi \in C^{\infty}(\Re;[0,1])$ be a smooth function with $\int_0^1 \varphi(s)ds > 0$, $\varphi(s) = 0$ for all $s \le 0$ and $s \ge 1$. Define $\eta(t) := \int_0^1 \varphi(s)\overline{\eta}(t+s)ds$, which is a $C^{\infty}$ positive function that satisfies $\eta(t) \le \tilde{\eta}(t)$ for all $t \ge 0$. Consequently, by virtue of (2.21) we obtain:

$$|x| \le 2\eta(t) \;\Rightarrow \Xi(t,x,0) \le 0 \tag{2.22}$$

Let $h \in C^{\infty}(\Re;[0,1])$ be a smooth non-decreasing function with $h(s) = 0$ for all $s \le 0$ and $h(s) = 1$ for all $s \ge 1$. We define:

$$k(t,x) := h\left(\frac{|x|^2 - 2\eta^2(t)}{2\eta^2(t)}\right)\sum_{i=1}^{\infty}\theta_i(t,x)u_i \tag{2.23}$$

Clearly, $k$ as defined by (2.23) is a smooth function with $k(t,0) = 0$ for all $t \ge 0$. Moreover, since $k(t,x)$ is defined as a (finite) convex combination of $u_i \in U$ and $0 \in U$, we have $k(t,x) \in U$ for all $(t,x) \in \Re^+ \times \Re^n$.

Let $(t,x) \in \Re^+ \times \Re^n$ with $|x| \ge 2\,\eta(t)$ and define $J(t,x) = \{ j \in \{1,2,...\}\,;\, \theta_j(t,x) \ne 0 \}$ (a finite set). Notice that by virtue of (2.9) and definition (2.23) we get:

$$\sup_{d\in D} V^0\big(t,x;\, f(t,d,x,k(t,x))\big) = \sup_{d\in D} V^0\left(t,x;\, f\left(t,d,x,\sum_{j\in J(t,x)}\theta_j(t,x)u_j\right)\right) \le -\rho(V(t,x)) + \max_{j\in J(t,x)}\{\Psi(t,x,u_j)\} \tag{2.24}$$

Notice that for each $j \in J(t,x)$ we obtain from (2.17) that $(t,x) \in \Omega_j$. Consequently, by virtue of (2.16) and definition (2.13a) we have that $\Psi(t,x,u_j) \le 8q(t)$, for all $j \in J(t,x)$. Combining the previous inequality with inequality (2.24), we conclude that the following property holds for all $(t,x,d) \in \Re^+ \times \Re^n \times D$ with $|x| \ge 2\eta(t)$:

$$V^0(t,x;\, f(t,d,x,k(t,x))) \le -\rho(V(t,x)) + 8q(t) \tag{2.25}$$

Let $(t,x) \in \Re^+ \times \Re^n$ with $|x| \le \sqrt{2}\,\eta(t)$. Notice that by virtue of definition (2.23) we get:

$$\sup_{d\in D} V^0\big(t,x;\, f(t,d,x,k(t,x))\big) = \sup_{d\in D} V^0\big(t,x;\, f(t,d,x,0)\big)$$

By virtue of (2.22), (2.13a) and the above inequality we conclude that (2.25) holds as well for all $(t,x) \in \Re^+ \times \Re^n$ with $|x| \le \sqrt{2}\,\eta(t)$. Finally, for the case $(t,x) \in \Re^+ \times \Re^n$ with $\sqrt{2}\,\eta(t) < |x| < 2\eta(t)$, let $J(t,x) = \{ j \in \{1,2,...\}\,;\, \theta_j(t,x) \ne 0 \}$ and notice that from (2.9) we get:

$$\sup_{d\in D} V^0\big(t,x;\, f(t,d,x,k(t,x))\big) =$$

$$\sup_{d\in D} V^0\left(t,x;\, f\left(t,d,x,h\left(\frac{|x|^2 - 2\eta^2(t)}{2\eta^2(t)}\right)\sum_{j\in J(t,x)}\theta_j(t,x)u_j\right)\right) \le -\rho(V(t,x)) + \max\{\Psi(t,x,0), \Psi(t,x,u_j),\, j \in J(t,x)\} \tag{2.26}$$

Taking into account definition (2.13a) and inequalities (2.16), (2.17), (2.22), (2.26) we may conclude that (2.25) holds as well for all $(t,x) \in \Re^+ \times \Re^n$ with $\sqrt{2}\;\eta(t) < |x| < 2\eta(t)$. Consequently, (2.25) holds for all $(t,x) \in \Re^+ \times \Re^n$.

It follows from (2.25) and Proposition 2.5 that system (1.1) with $u = k(t,x)$ is RGAOS.

**(b)** $\Rightarrow$ **(d)** Since system (1.1) with $u = k(t,x)$ is RGAOS, and since $f(t,d,x,k(t,x))$ is Lipschitz with respect to $x$ on each bounded subset of $\Re^+ \times \Re^n \times D$, it follows from Theorem 3.2 in [12] that there exists a function $V \in C^\infty(\Re^+ \times \Re^n; \Re^+)$, functions $a_1, a_2 \in K_\infty$, $\beta, \mu \in K^+$ such that

$$a_1\left(\left|(\mu(t)x, H(t,x))\right|\right) \le V(t,x) \le a_2\left(\beta(t)|x|\right),\ \forall (t,x) \in \Re^+ \times \Re^n \tag{2.27a}$$

$$\frac{\partial V}{\partial t}(t,x) + \sup_{d \in D}\left( \frac{\partial V}{\partial x}(t,x) f(t,d,x,k(t,x)) \right) \le -V(t,x),\ \forall (t,x) \in \Re^+ \times \Re^n \tag{2.27b}$$

We next prove that $V$ is an ORCLF for (1.1). Obviously property (i) of Definition 2.6 is a consequence of inequality (2.27a). Define

$$\Psi(t,x,u) := V(t,x) + \sup\left\{ \frac{\partial V}{\partial t}(t,x) + \frac{\partial V}{\partial x}(t,x) f(t,d,x,v) : d \in D,\ v \in U \ with\ |v - k(t,x)| \le |u - k(t,x)| \right\} \tag{2.28}$$

Inequality (2.8) with $q(t) \equiv 0$ is an immediate consequence of inequality (2.27b) and definition (2.28). Moreover, it holds that

$$\frac{\partial V}{\partial t}(t,x) + \sup_{d \in D}\left( \frac{\partial V}{\partial x}(t,x) f(t,d,x,u) \right) \le -V(t,x) + \Psi(t,x,u),\ \forall (t,x,u) \in \Re^+ \times \Re^n \times U \tag{2.29}$$

It follows from compactness of $D \subset \Re^l$, continuity of $V(t,x)$, $\frac{\partial V}{\partial t}(t,x) + \frac{\partial V}{\partial x}(t,x) f(t,d,x,v)$, $k : \Re^+ \times \Re^n \to U$ and Theorem 1.4.16 in [2] that the function $\Psi$ as defined by (2.28) is upper semi-continuous. Clearly, definition (2.28) implies $\Psi(t,0,0) = 0$ for all $t \ge 0$.

Finally, we show that inequality (2.9) holds with $\rho(s) := s$. Let arbitrary $(t,x) \in \Re^+ \times \Re^n$, $\{u_1, u_2, ..., u_p\} \subset U$ and $\lambda_i \in [0,1]$ ($i = 1,...,p$) with $\sum_{i=1}^{p} \lambda_i = 1$. Definition (2.28) implies:

$$\begin{aligned}
&\Psi\left(t,x,\sum_{i=1}^{p}\lambda_i u_i\right) \\
&= V(t,x) + \sup\left\{ \frac{\partial V}{\partial t}(t,x) + \frac{\partial V}{\partial x}(t,x) f(t,d,x,v) : d \in D,\ v \in U \ with\ |v - k(t,x)| \le \left|\sum_{i=1}^{p}\lambda_i u_i - k(t,x)\right| \right\} \\
&\le V(t,x) + \sup\left\{ \frac{\partial V}{\partial t}(t,x) + \frac{\partial V}{\partial x}(t,x) f(t,d,x,v) : d \in D,\ v \in U \ with\ |v - k(t,x)| \le \sum_{i=1}^{p}\lambda_i |u_i - k(t,x)| \right\} \\
&\le V(t,x) + \sup\left\{ \frac{\partial V}{\partial t}(t,x) + \frac{\partial V}{\partial x}(t,x) f(t,d,x,v) : d \in D,\ v \in U \ with\ |v - k(t,x)| \le \max_{i=1,...,p} |u_i - k(t,x)| \right\}
\end{aligned}$$

Let $j \in \{1,...,p\}$ such that $|u_j - k(t,x)| = \max_{i=1,...,p} |u_i - k(t,x)|$. The previous inequalities imply that

$$\Psi\left(t,x,\sum_{i=1}^{p}\lambda_i u_i\right) \le \Psi(t,x,u_j) \le \max_{i=1,...,p} \Psi(t,x,u_i) \tag{2.30}$$

Inequality (2.9) with $\rho(s) := s$ is a direct consequence of (2.29) and (2.30).

**(a)** $\Rightarrow$ **(c)** The proof is exactly the same with the proof of implication (b) $\Rightarrow$ (d). The only additional point is that by virtue of Lemma 3.2 in [11] there exist functions $a \in K_\infty$, $\gamma \in K^+$ such that

$$|k(t,x)| \le a\left(\gamma(t)|x|\right),\ \forall (t,x) \in \Re^+ \times \Re^n \tag{2.31}$$

Consequently, for every $(t,x)\in \Re^{+}\times \Re^{n}$ there exists $u\in U$ with $|u|\le a\big(\gamma(t)|x|\big)$ (namely $u=k(t,x)$) such that (2.10) holds with $\rho(s):=s$ and $q(t)\equiv 0$. The proof is complete. $\lhd$

**Proof of Theorem 2.9:** Suppose that (1.1) admits an ORCLF which satisfies the small control property with $q(t)\equiv 0$. Define:

$$\Xi(t,x,u):=\Psi(t,x,u)-\frac{1}{2}\rho\big(V(t,x)\big),\ (t,x,u)\in \Re^{+}\times \Re^{n}\times U \tag{2.32a}$$

$$\Xi(t,x,u):=\Xi(0,x,u),\ (t,x,u)\in(-1,0)\times \Re^{n}\times U \tag{2.32b}$$

The definition of $\Xi$, given by (2.32a,b), guarantees that the function $\Xi:(-1,+\infty)\times\Re^{n}\times U\to\Re\cup\{+\infty\}$ with $\Xi(t,0,0)=0$ for all $t>-1$ is such that, for each $u\in U$ the mapping $(t,x)\to\Xi(t,x,u)$ is upper semi-continuous. By virtue of (2.10) with $q(t)\equiv 0$ and upper semi-continuity of $\Xi$, it follows that for each $(t,x)\in(-1,+\infty)\times(\Re^{n}\setminus\{0\})$ there exist $u=u(t,x)\in U$ with $|u|\le a\big(\gamma(\max(0,t))|x|\big)$ and $\delta=\delta(t,x)\in(0,\min(1,t+1))$, $\delta(t,x)\le\frac{|x|}{2}$, such that

$$\Xi(\tau,y,u(t,x))\le 0,\ \forall(\tau,y)\in\big\{(\tau,y)\in(-1,+\infty)\times\Re^{n}:|\tau-t|+|y-x|<\delta\big\} \tag{2.33}$$

Using (2.33) and standard partition of unity arguments, we can determine sequences $\{(t_i,x_i)\in(-1,+\infty)\times(\Re^{n}\setminus\{0\})\}_{i=1}^{\infty}$, $\{u_i\in U\}_{i=1}^{\infty}$, $\{\delta_i\}_{i=1}^{\infty}$ with $|u_i|\le a\big(\gamma(\max(0,t_i))|x_i|\big)$, $\delta_i=\delta(t_i,x_i)\in(0,\min(1,t_i+1))$, $\delta_i=\delta(t_i,x_i)\le\frac{|x_i|}{2}$ associated with a sequence of open sets $\{\Omega_i\}_{i=1}^{\infty}$ with

$$\Omega_i\subseteq\big\{(\tau,y)\in(-1,+\infty)\times\Re^{n}:|\tau-t_i|+|y-x_i|<\delta_i\big\} \tag{2.34}$$

forming a locally finite open covering of $(-1,+\infty)\times(\Re^{n}\setminus\{0\})$ and in such a way that:

$$\Xi(\tau,y,u_i)\le 0,\ \forall(\tau,y)\in\Omega_i \tag{2.35}$$

Also, a family of smooth functions $\{\theta_i\}_{i=1}^{\infty}$ with $\theta_i(t,x)\ge 0$ for all $(t,x)\in(-1,+\infty)\times(\Re^{n}\setminus\{0\})$ can be determined with

$$s\,upp\ \theta_i\subseteq\Omega_i \tag{2.36}$$

$$\sum_{i=1}^{\infty}\theta_i(t,x)=1,\ \forall(t,x)\in(-1,+\infty)\times(\Re^{n}\setminus\{0\}) \tag{2.37}$$

We define:

$$k(t,x):=\sum_{i=1}^{\infty}\theta_i(t,x)u_i\ \text{ for } t\ge 0,\ x\ne 0 \tag{2.38a}$$

$$k(t,0):=0\ \text{ for } t\ge 0 \tag{2.38b}$$

It follows from (2.38a) that $k$ is $C^{\infty}$ on the set $\Re^{+}\times(\Re^{n}\setminus\{0\})$. Moreover, since $k(t,x)$ is defined as a (finite) convex combination of $u_i\in U$ and $0\in U$, we have $k(t,x)\in U$ for all $(t,x)\in\Re^{+}\times\Re^{n}$. In order to prove continuity of $k$ at zero, let $(t,x)\in\Re^{+}\times(\Re^{n}\setminus\{0\})$ and define $J(t,x)=\big\{j\in\{1,2,...\}\,;\theta_j(t,x)\ne 0\big\}$ (a finite set). Notice that for each $j\in J(t,x)$ we obtain from (2.36) that $(t,x)\in\Omega_j$. Consequently, using (2.34), (2.38a) and the facts $|u_j|\le a\big(\gamma(\max(0,t_j))|x_j|\big)$, $\delta_j=\delta(t_j,x_j)\in(0,\min(1,t_j+1))$, $\delta_j=\delta(t_j,x_j)\le\frac{|x_j|}{2}$, we obtain:

$$\left|k(t,x)\right| \leq \max_{j\in J(t,x)} \left|u_j\right| \leq \max_{j\in J(t,x)} a\left(\gamma(\max(0,t_j))\left|x_j\right|\right) \leq \tilde{a}\left(\tilde{\gamma}(t)\left|x\right|\right)$$

where $\tilde{a}(s):=a(2s)$ and $\tilde{\gamma}(t):=\max_{0\leq\tau\leq t+1}\gamma(\tau)$. The above inequality in conjunction with definition (2.38b) shows continuity of $k$ at zero. Next we show that:

$$V^0(t,x;\, f(t,d,x,k(t,x))) \leq -\frac{1}{2}\rho\left(V(t,x)\right),\ \forall (t,x,d)\in \Re^+\times\Re^n\times D \tag{2.39}$$

Clearly, by virtue of definition (2.38b) and inequality (2.3), it follows that inequality (2.39) holds for all $t\geq 0$, $x=0$. For $(t,x)\in\Re^+\times(\Re^n\setminus\{0\})$, define $J(t,x)=\left\{j\in\{1,2,...\}\,;\,\theta_j(t,x)\neq 0\right\}$ (a finite set). Notice that by virtue of (2.9) and definition (2.38a) we get:

$$\sup_{d\in D} V^0\left(t,x;\, f\left(t,d,x,k(t,x)\right)\right) = \sup_{d\in D} V^0\left(t,x;\, f\left(t,d,x,\sum_{j\in J(t,x)}\theta_j(t,x)u_j\right)\right) \leq -\rho\left(V(t,x)\right) + \max_{j\in J(t,x)}\left\{\Psi(t,x,u_j)\right\} \tag{2.40}$$

Notice that for each $j\in J(t,x)$ we obtain from (2.36) that $(t,x)\in\Omega_j$. Consequently, by virtue of (2.35) and definition (2.32a) we have that $\Psi(t,x,u_j)\leq\frac{1}{2}\rho\left(V(t,x)\right)$, for all $j\in J(t,x)$. Combining the previous inequality with inequality (2.40), we conclude that (2.39) holds.

If the ORCLF $V$ and the function $\Psi$ involved in property (ii) of Definition 2.6 are time independent then the partition of unity arguments used above may be repeated on $\Re^n\setminus\{0\}$ instead of $\Re^+\times(\Re^n\setminus\{0\})$. This implies that the constructed feedback is time invariant.

In order to show uniqueness of solutions for the closed-loop system (2.11) we consider the dynamical system

$$\dot{x}(t) = f(t,d(t),x(t),k(t,x(t))) \quad , \quad x(t)\in\Re^n\setminus\{0\} \tag{2.41}$$

It is clear from hypothesis (H1) and smoothness of $k$ on the set $\Re^+\times(\Re^n\setminus\{0\})$, that for every $(t_0,x_0,d)\in\Re^+\times(\Re^n\setminus\{0\})\times M_D$, the solution with initial condition $x(t_0)=x_0\in\Re^n\setminus\{0\}$, corresponding to $d\in M_D$ is unique and is defined on the interval $[t_0,t_{\max})$, where $t_{\max}>t_0$ is the maximal existence time of the solution of (2.41).

Notice that the solution of (2.11) with initial condition $x(t_0)=x_0\in\Re^n\setminus\{0\}$, corresponding to some $d\in M_D$ coincides with the unique solution of (2.41) evolving on $\Re^+\times(\Re^n\setminus\{0\})$ with same initial condition $x(t_0)=x_0\in\Re^n\setminus\{0\}$, and same $d\in M_D$ on the interval $[t_0,t_{\max})$, where $t_{\max}>t_0$ is the maximal existence time of the solution of (2.41).

For the case $t_{\max}=+\infty$, uniqueness of solutions for (2.11) is a direct consequence of previous argument. Suppose next that $t_{\max}<+\infty$. To establish uniqueness of solutions for (2.11), we need the following implication, which is a consequence of (2.6) and (2.39):

$$t_{\max}<+\infty \Rightarrow \lim_{t\to t_{\max}^-} x(t)=0 \tag{2.42}$$

In order to show (2.42), let $(t_0,x_0,d)\in\Re^+\times(\Re^n\setminus\{0\})\times M_D$ and suppose that the maximal existence time $t_{\max}>t_0$ of the (unique) solution of (2.41) with initial condition $x(t_0)=x_0\in\Re^n\setminus\{0\}$ corresponding to $d\in M_D$ is finite, i.e., $t_{\max}<+\infty$. Lemma 2.5 in conjunction with (2.39) implies that

$$V(t,x(t))\leq V(t_0,x_0)\,,\ \forall t\in[t_0,t_{\max}) \tag{2.43}$$

The above inequality in conjunction with (2.6) with $\beta(t)\equiv 1$ gives

$$|x(t)| \le M := \frac{1}{\min\limits_{\tau \in [0,t_{\max}]} \mu(\tau)} a_1^{-1}\left(a_2(|x_0|)\right) < +\infty \, , \ \forall t \in [t_0, t_{\max}) \tag{2.44}$$

Definition of $t_{\max}$ and (2.44) implies (2.42). By applying standard arguments we may also establish show that for every $(t_0, d) \in \Re^+ \times M_D$, the solution of (2.11) with initial condition $x(t_0) = 0$, corresponding to input $d \in M_D$ is unique and satisfies $x(t) = 0$ for all $t \ge t_0$. Indeed, this follows from Lemma 2.5 and inequality (2.39), which imply inequality (2.43). The previous discussion in conjunction with (2.42) asserts that the solution $x(\cdot)$ of (2.11) with initial condition $x(t_0) = x_0 \in \Re^n \setminus \{0\}$, corresponding to $d \in M_D$ coincides with the solution of (2.41) with same initial condition, and same $d \in M_D$ on the interval $[t_0, t_{\max})$, $t_{\max} > t_0$ being the maximal existence time of the solution (2.41); moreover, if $t_{\max} < +\infty$, the corresponding solution of (2.11) satisfies $x(t) = 0$ for all $t \ge t_{\max}$ and uniqueness of solutions for (2.11) is established.

The fact that (2.11) is URGAOS follows directly from Proposition 2.5 and inequality (2.39).

Finally, if there exist functions $\eta \in K^+$, $\varphi \in C^v(A;U)$ where $A = \bigcup\limits_{t \ge 0} \{t\} \times B_t$, $B_t = \{x \in \Re^n : |x| < 4\eta(t)\}$ with $\varphi(t,0) = 0$ for all $t \ge 0$, such that (2.12) holds, then we consider the smooth feedback defined by:

$$k(t,x) := \left(1 - h\left(\frac{|x|^2 - 2\eta^2(t)}{2\eta^2(t)}\right)\right)\varphi\left(t, \Pr_{Q(t)}(x)\right) + h\left(\frac{|x|^2 - 2\eta^2(t)}{2\eta^2(t)}\right)\tilde{k}(t,x) \tag{2.45}$$

where

$$\tilde{k}(t,x) := \sum_{i=1}^{\infty} \theta_i(t,x) u_i \ \text{ for } t \ge 0, \ x \ne 0 \tag{2.46a}$$

$$\tilde{k}(t,0) := 0 \ \text{ for } t \ge 0 \tag{2.46b}$$

$h \in C^\infty(\Re;[0,1])$ is a smooth non-decreasing function with $h(s) = 0$ for all $s \le 0$ and $h(s) = 1$ for all $s \ge 1$ and $Q(t) = \{x \in \Re^n : |x| \le 3\eta(t)\}$. Clearly, $k$ as defined by (2.45) is of class $C^v(\Re^+ \times \Re^n; U)$ with $k(t,0) = 0$ for all $t \ge 0$. Using the same arguments as in the proof of Theorem 2.8 we may establish that (2.39) holds. The proof is complete. $\lhd$

## 3. Additional Remarks and Examples on the Finite-Dimensional Case

The problem with Definition 2.1 of the ORCLF that might arise in practice is the assumption of the knowledge of the function $\Psi : \Re^+ \times \Re^n \times U \to \Re \cup \{+\infty\}$ involved in property (ii) of Definition 2.6. Particularly, the following problem arises:

**Problem (P):** *Consider system (1.1) under hypotheses (H1-3) and the following hypothesis:*

**(H4)** *There exist a $C^1$ function $V : \Re^+ \times \Re^n \to \Re^+$, which satisfies property (i) of Definition 2.6, a function $q \in \mathfrak{E}$ and a $C^0$ positive definite function $\rho : \Re^+ \to \Re^+$ such that the following inequality holds:*

$$\inf_{u \in U} \left\{ \frac{\partial V}{\partial t}(t,x) + \sup_{d \in D}\left( \frac{\partial V}{\partial x}(t,x) f(t,d,x,u) \right)\right\} \le -\rho\left(V(t,x)\right) + q(t), \ \forall (t,x) \in \Re^+ \times \Re^n \tag{3.1}$$

*Is $V : \Re^+ \times \Re^n \to \Re^+$ an ORCLF for system (1.1)?*

The proof of implication (b)$\Rightarrow$(d) of Theorem 2.8 gives the solution to Problem (P): If there exists a continuous function $k : \Re^+ \times \Re^n \to U$ with $k(t,0) = 0$ for all $t \ge 0$, a function $\tilde{q} \in \mathfrak{E}$ and a $C^0$ positive definite function

$\tilde{\rho}:\Re^+ \to \Re^+$ such that

$$\frac{\partial V}{\partial t}(t,x)+\sup_{d\in D}\left(\frac{\partial V}{\partial x}(t,x)\,f(t,d,x,k(t,x))\right)\le -\tilde{\rho}\big(V(t,x)\big)+\tilde{q}(t)\,,\ \forall (t,x)\in\Re^+\times\Re^n \tag{3.2}$$

then $V$ is an ORCLF for (1.1). Particularly, the function $\Psi:\Re^+\times\Re^n\times U\to\Re\cup\{+\infty\}$ involved in property (ii) of Definition 2.6 may be defined by

$$\Psi(t,x,u):=\tilde{\rho}\big(V(t,x)\big)+\sup\left\{\frac{\partial V}{\partial t}(t,x)+\frac{\partial V}{\partial x}(t,x)\,f(t,d,x,v): d\in D\,,\, v\in U\ with\ |v-k(t,x)|\le|u-k(t,x)|\right\} \tag{3.3}$$

The reader may check that $\Psi:\Re^+\times\Re^n\times U\to\Re\cup\{+\infty\}$ as defined by (3.3) satisfies inequalities (2.8), (2.9) of Definition 2.6, following exactly the same procedure as in the proof of implication (b)$\Rightarrow$(d) of Theorem 2.8. Moreover, by virtue of Theorem 2.8, if $V:\Re^+\times\Re^n\to\Re^+$ an ORCLF for system (1.1) then the proof of implication (d)$\Rightarrow$(a) of Theorem 2.8 shows that there exists a continuous function $k:\Re^+\times\Re^n\to U$ with $k(t,0)=0$ for all $t\ge 0$, a function $\tilde{q}\in\mathcal{E}$ and a $C^0$ positive definite function $\tilde{\rho}:\Re^+\to\Re^+$ such that (3.2) holds. Consequently, $V:\Re^+\times\Re^n\to\Re^+$ is an ORCLF for (1.1) under hypotheses (H1-4) **if and only if** there exists a continuous function $k:\Re^+\times\Re^n\to U$ with $k(t,0)=0$ for all $t\ge 0$, a function $\tilde{q}\in\mathcal{E}$ and a $C^0$ positive definite function $\tilde{\rho}:\Re^+\to\Re^+$ such that (3.2) holds.

The problem with the above solution to Problem (P) is that we can check if $V:\Re^+\times\Re^n\to\Re^+$ is an ORCLF for (1.1) by constructing a feedback stabilizer for (1.1). On the other hand, our goal in practice is to construct the feedback stabilizer based on the mere knowledge of the Lyapunov function $V:\Re^+\times\Re^n\to\Re^+$ under hypotheses (H1-4). Consequently, the above solution to Problem (P) is useless for feedback construction purposes.

The rest of the section provides sufficient conditions for establishing that $V:\Re^+\times\Re^n\to\Re^+$ under hypotheses (H1-4) is an ORCLF for (1.1).

Indeed, if the mapping $u\to\frac{\partial V}{\partial x}(t,x)\,f(t,d,x,u)$ is quasi-convex for each fixed $(t,x,d)\in\Re^+\times\Re^n\times D$ then the mapping $u\to\frac{\partial V}{\partial t}(t,x)+\sup_{d\in D}\left(\frac{\partial V}{\partial x}(t,x)\,f(t,d,x,u)\right)$ is quasi-convex for each fixed $(t,x)\in\Re^+\times\Re^n$. It follows that property (ii) of Definition 2.6 is satisfied with $\Psi(t,x,u):=\rho\big(V(t,x)\big)+\frac{\partial V}{\partial t}(t,x)+\sup_{d\in D}\left(\frac{\partial V}{\partial x}(t,x)\,f(t,d,x,u)\right)$.
This is exactly the case arising in affine in the control systems: for affine in the control systems the mapping $u\to\frac{\partial V}{\partial x}(t,x)\,f(t,d,x,u)$ is convex.

The following lemma helps us to generalize the above sufficient condition.

**Lemma 3.1:** *Let the mapping* $f:A\times U\to\Re$, *where* $U\subseteq\Re^m$ *a closed convex set. Define the set-valued map:*

$$A\times U\ni(x,u)\to\mathcal{U}(x,u):=\overline{co}\{\,v\in U:\ f(x,v)\le f(x,u)\,\} \tag{3.4}$$

*and the mapping* $\psi:A\times U\to\Re\cup\{+\infty\}$

$$\psi(x,u):=\sup\{\,f(x,v):v\in\mathcal{U}(x,u)\,\} \tag{3.5}$$

*Then for every finite set* $\{u_1,u_2,...,u_p\}\subset U$ *and for every* $\lambda_i\in[0,1]$ *(* $i=1,...,p$ *) with* $\sum_{i=1}^{p}\lambda_i=1$, *it holds that:*

$$f\left(x,\sum_{i=1}^{p}\lambda_i u_i\right)\le\max\{\psi(x,u_i)\,,\,i=1,...,p\,\},\ \forall x\in A \tag{3.6}$$

**Proof:** Let a finite set $\{u_1, u_2, ..., u_p\} \subset U$ and $\lambda_i \in [0,1]$ ($i = 1,..., p$) with $\sum_{i=1}^{p} \lambda_i = 1$. Let $u \in \{u_1, u_2, ..., u_p\}$ with $f(x,u) = \max_{i=1,...,p} f(x,u_i)$. It follows from definition (3.4) that $\sum_{i=1}^{p} \lambda_i u_i \in \mathcal{U}(x,u)$ and consequently, by virtue of definition (3.5), we get $f\left(x, \sum_{i=1}^{p} \lambda_i u_i\right) \le \psi(x,u)$. The previous inequality combined with the fact that $u \in \{u_1, u_2, ..., u_p\}$ (which implies $\psi(x,u) \le \max\{\psi(x,u_i), i = 1,..., p\}$) establishes (3.6). The proof is complete. $\lhd$

The following lemma is a direct consequence of the previous lemma.

**Lemma 3.2:** *Let $V : \Re^+ \times \Re^n \to \Re^+$ be a $C^1$ function which satisfies the following properties:*

**(i)** *there exists a function $q \in \mathcal{E}$ and a $C^0$ positive definite function $\rho : \Re^+ \to \Re^+$ such that for every $(t,x) \in \Re^+ \times \Re^n$ there exists $u \in U$ with $\mathcal{U}(t,x,u) \subseteq \tilde{\mathcal{U}}(t,x) \neq \emptyset$, where*
$$\mathcal{U}(t,x,u) := \overline{co}\left\{ v \in U : \sup_{d\in D}\left( \frac{\partial V}{\partial x}(t,x) f(t,d,x,v) \right) \le \sup_{d\in D}\left( \frac{\partial V}{\partial x}(t,x) f(t,d,x,u) \right) \right\}$$
*and*
$$\tilde{\mathcal{U}}(t,x) := \left\{ v \in U : \frac{\partial V}{\partial t}(t,x) + \sup_{d\in D}\left( \frac{\partial V}{\partial x}(t,x) f(t,d,x,v) \right) \le -\rho(V(t,x)) + q(t) \right\},$$

**(ii)** *for each fixed $u \in U$ the mapping $\Re^+ \times \Re^n \ni (t,x) \to \sup\left\{ \sup_{d\in D}\left( \frac{\partial V}{\partial x}(t,x) f(t,d,x,v) \right) : v \in \mathcal{U}(t,x,u) \right\}$ is upper semi-continuous.*

*Then property (ii) of Definition 2.6 holds with*
$$\Psi(t,x,u) := \rho(V(t,x)) + \frac{\partial V}{\partial t}(t,x) + \sup\left\{ \sup_{d\in D}\left( \frac{\partial V}{\partial x}(t,x) f(t,d,x,v) \right) : v \in \mathcal{U}(t,x,u) \right\}.$$

It should be noted that if the mapping $u \to \frac{\partial V}{\partial x}(t,x) f(t,d,x,u)$ is quasi-convex for each fixed $(t,x,d) \in \Re^+ \times \Re^n \times D$ then the set-valued map $\mathcal{U}(t,x,u) := \overline{co}\left\{ v \in U : \sup_{d\in D}\left( \frac{\partial V}{\partial x}(t,x) f(t,d,x,v) \right) \le \sup_{d\in D}\left( \frac{\partial V}{\partial x}(t,x) f(t,d,x,u) \right) \right\}$ in property (i) of Lemma 3.2 satisfies $\mathcal{U}(t,x,u) = \left\{ v \in U : \sup_{d\in D}\left( \frac{\partial V}{\partial x}(t,x) f(t,d,x,v) \right) \le \sup_{d\in D}\left( \frac{\partial V}{\partial x}(t,x) f(t,d,x,u) \right) \right\}$ and consequently property (i) of Lemma 3.2 becomes equivalent to the existence of $u \in U$ with $\frac{\partial V}{\partial t}(t,x) + \sup_{d\in D}\left( \frac{\partial V}{\partial x}(t,x) f(t,d,x,u) \right) \le -\rho(V(t,x)) + q(t)$.

The following example illustrates the use of Lemma 3.2 for a special class of nonlinear systems.

**Example 3.3:** Consider system (1.1) under hypotheses (H1-3) with $m = 1$, $U = \Re$ and a $C^1$ function $V : \Re^+ \times \Re^n \to \Re^+$ which satisfies property (i) of Definition 2.6 as well as

$$\frac{\partial V}{\partial t}(t,x) + \sup_{d\in D}\left( \frac{\partial V}{\partial x}(t,x) f(t,d,x,u) \right) = a(t,x)u^2 + b(t,x)u + c(t,x), \ \forall (t,x,u) \in \Re^+ \times \Re^n \times U \tag{3.7}$$

$$\inf_{u\in\Re}\left(a(t,x)u^2+b(t,x)u+c(t,x)\right)\le-\rho\left(V(t,x)\right)+q(t)\,,\ \forall(t,x)\in\Re^+\times\Re^n \tag{3.8}$$

for appropriate continuous mappings $a,b,c:\Re^+\times\Re^n\to\Re$ with $a(t,0)=b(t,0)=c(t,0)=0$ for all $t\ge 0$, a function $q\in\mathcal{E}$ and a $C^0$ positive definite function $\rho:\Re^+\to\Re^+$. We next prove that $V:\Re^+\times\Re^n\to\Re^+$ is an ORCLF for (1.1) provided that the following implications hold:

$$a(t,x)<0\quad\Rightarrow\quad -\frac{b^2(t,x)}{4a(t,x)}+c(t,x)\le-\rho\left(V(t,x)\right)+q(t) \tag{3.9}$$

for every sequence $\{(t_i,x_i)\}_{i=0}^{\infty}$ with $a(t_i,x_i)<0$ and $(t_i,x_i)\to(t,x)$ with $a(t,x)=0$ it holds that $\dfrac{b^2(t_i,x_i)}{|a(t_i,x_i)|}\to 0$

(3.10)

Following the notation of Lemma 3.2 we define $\mathcal{U}(t,x,u)=\overline{co}\left\{v\in\Re: a(t,x)v^2+b(t,x)v\le a(t,x)u^2+b(t,x)u\right\}$ and $\tilde{\mathcal{U}}(t,x)=\left\{v\in\Re: a(t,x)v^2+b(t,x)v+c(t,x)\le-\rho(V(t,x))+q(t)\right\}$. We notice that:

a) If $a(t,x)>0$ then for every $u\in\Re$ the set $\left\{v\in\Re: a(t,x)v^2+b(t,x)v\le a(t,x)u^2+b(t,x)u\right\}$ is closed and convex and consequently, there exists $u\in\Re$ with $\mathcal{U}(t,x,u)\subseteq\tilde{\mathcal{U}}(t,x)$ (specifically, by virtue of (3.8), the inclusion $\mathcal{U}(t,x,u)\subseteq\tilde{\mathcal{U}}(t,x)$ holds for $u=-\dfrac{b(t,x)}{2a(t,x)}$).

b) If $a(t,x)=0$ then for every $u\in\Re$ the set $\left\{v\in\Re: b(t,x)v\le b(t,x)u\right\}$ is closed and convex and consequently, there exists $u\in\Re$ with $\mathcal{U}(t,x,u)\subseteq\tilde{\mathcal{U}}(t,x)$. Specifically, if $b(t,x)\ne 0$ then the inclusion $\mathcal{U}(t,x,u)\subseteq\tilde{\mathcal{U}}(t,x)$ holds for $u=-\dfrac{c(t,x)+\rho(V(t,x))}{b(t,x)}$. If $b(t,x)=0$, then by virtue of (3.8), the inclusion $\mathcal{U}(t,x,u)\subseteq\tilde{\mathcal{U}}(t,x)$ holds for every $u\in\Re$.

c) If $a(t,x)<0$ then for every $u\ne-\dfrac{b(t,x)}{2a(t,x)}$ the set $\left\{v\in\Re: a(t,x)v^2+b(t,x)v\le a(t,x)u^2+b(t,x)u\right\}$ is not convex and there exist $v_1,v_2\in\Re$ with $v_1<v_2$ such that $\left\{v\in\Re: a(t,x)v^2+b(t,x)v\le a(t,x)u^2+b(t,x)u\right\}=(-\infty,v_1]\cup[v_2,+\infty)$. On the other hand if $u=-\dfrac{b(t,x)}{2a(t,x)}$ then it holds that $\left\{v\in\Re: a(t,x)v^2+b(t,x)v\le a(t,x)u^2+b(t,x)u\right\}=\Re$. Consequently, if $a(t,x)<0$ then for every $u\in\Re$ it holds that $\mathcal{U}(t,x,u)=\Re$. However, in this case implication (3.9) guarantees that $\tilde{\mathcal{U}}(t,x)=\Re$ and therefore the inclusion $\mathcal{U}(t,x,u)\subseteq\tilde{\mathcal{U}}(t,x)$ holds for every $u\in\Re$.

Thus, property (i) of Lemma 3.2 holds for the function $V:\Re^+\times\Re^n\to\Re^+$. Since $\Psi(t,x,u):=\rho\left(V(t,x)\right)+\dfrac{\partial V}{\partial t}(t,x)+\sup\left\{\sup_{d\in D}\left(\dfrac{\partial V}{\partial x}(t,x)f(t,d,x,v)\right): v\in\mathcal{U}(t,x,u)\right\}$, by virtue of all the above specifications for the set-valued map $\mathcal{U}(t,x,u)$ we get:

$$\Psi(t,x,u):=\rho\left(V(t,x)\right)+\begin{cases}a(t,x)u^2+b(t,x)u+c(t,x) & if\quad a(t,x)\ge 0\\ -\dfrac{b^2(t,x)}{4a(t,x)}+c(t,x) & if\quad a(t,x)<0\end{cases} \tag{3.11}$$

Notice that implication (3.10) guarantees that property (ii) of Lemma 3.2 holds and consequently property (ii) of Definition 2.6 holds with $\Psi$ defined by (3.11).
The reader should notice that other choices for the mapping $\Psi(t,x,u)$ are possible. For example, the selection

$$\Psi(t,x,u):=\rho\left(V(t,x)\right)+\begin{cases}a(t,x)u^2+b(t,x)u+c(t,x) & if\quad a(t,x)\ge 0\\ b(t,x)u+c(t,x) & if\quad a(t,x)<0\end{cases} \tag{3.12}$$

guarantees that $V:\Re^+ \times \Re^n \to \Re^+$ is an ORCLF for (1.1) provided that (3.7), (3.8) as well as the following implication holds:

$$a(t,x)<0\,,\,b(t,x)=0 \quad \Rightarrow \quad c(t,x) \le -\rho\big(V(t,x)\big)+q(t) \tag{3.13}$$

Notice that if (3.13) holds then property (ii) of Definition 2.6 holds with $\Psi$ defined by (3.12). Moreover, notice that if implication (3.9) holds then implication (3.13) automatically holds. $\lhd$

The following lemma provides a "patchy" construction by combining the formula provided by Lemma 3.2 and the knowledge of appropriate functions that can be used in certain regions of $\Re^+ \times \Re^n$ as feedback functions.

**Lemma 3.4:** *Let $V:\Re^+ \times \Re^n \to \Re^+$ be a $C^1$ function and suppose that there exist sets $\Omega_i \subseteq \Re^+ \times \Re^n$ ($i=0,...,p$) with $\Omega_i \cap \Omega_j = \emptyset$ for $i \ne j$ and $\bigcup_{i=0,...,p} \Omega_i = \Re^+ \times \Re^n$, functions $k_i:\Omega_i \to U$ ($i=1,...,p$), a function $q \in \mathcal{E}$ and a $C^0$ positive definite function $\rho:\Re^+ \to \Re^+$ such that:*

**(i)** *for every $(t,x)\in\Omega_0$ there exists $u \in U$ with $\mathcal{U}(t,x,u) \subseteq \tilde{\mathcal{U}}(t,x) \ne \emptyset$, where*

$$\mathcal{U}(t,x,u) := \overline{co}\left\{ v \in U : \sup_{d\in D}\left( \frac{\partial V}{\partial x}(t,x) f(t,d,x,v) \right) \le \sup_{d\in D}\left( \frac{\partial V}{\partial x}(t,x) f(t,d,x,u) \right) \right\} \quad \textit{and}$$

$$\tilde{\mathcal{U}}(t,x) := \left\{ v \in U : \frac{\partial V}{\partial t}(t,x) + \sup_{d\in D}\left( \frac{\partial V}{\partial x}(t,x) f(t,d,x,v) \right) \le -\rho\big(V(t,x)\big)+q(t) \right\},$$

**(ii)** *for every $i=1,...,p$ and $(t,x)\in\Omega_i$ it holds that*

$$\frac{\partial V}{\partial t}(t,x) + \sup_{d\in D}\left( \frac{\partial V}{\partial x}(t,x) f(t,d,x,k_i(t,x)) \right) \le -\rho\big(V(t,x)\big)+q(t)$$

*Consider the function $\Psi:\Re^+ \times \Re^n \times U \to \Re \cup \{+\infty\}$ defined by:*

$$\Psi(t,x,u) := \rho\big(V(t,x)\big) + \frac{\partial V}{\partial t}(t,x) + \sup\left\{ \sup_{d\in D}\left( \frac{\partial V}{\partial x}(t,x) f(t,d,x,v) \right) : v \in \mathcal{U}(t,x,u) \right\}, \text{ for } (t,x)\in\Omega_0 \tag{3.14a}$$

$$\Psi(t,x,u) := \rho\big(V(t,x)\big) + \sup\left\{ \frac{\partial V}{\partial t}(t,x) + \frac{\partial V}{\partial x}(t,x) f(t,d,x,v) : d \in D\,,\, v \in U \textit{ with } |v-k_i(t,x)| \le |u-k_i(t,x)| \right\},$$
$$\text{for } (t,x)\in\Omega_i\,,\ i=1,...,p \tag{3.14b}$$

*and suppose that $\Psi:\Re^+ \times \Re^n \times U \to \Re \cup \{+\infty\}$ is upper semi-continuous.*

*Then property (ii) of Definition 2.6 holds with $\Psi:\Re^+ \times \Re^n \times U \to \Re \cup \{+\infty\}$ as defined by (3.14) and $V:\Re^+ \times \Re^n \to \Re^+$ is an ORCLF for (1.1).*

The following example illustrates the efficiency of Lemma 3.4. It shows that the knowledge of appropriate functions that can be used in certain regions of $\Re^+ \times \Re^n$ as feedback functions, helps us to obtain less conservative results.

**Example 3.5:** Consider system (1.1) under hypotheses (H1-3) with $m=1$, $U=\Re$ and a $C^1$ function $V:\Re^+ \times \Re^n \to \Re^+$ which satisfies property (i) of Definition 2.6 as well as (3.7), (3.8) for appropriate continuous mappings $a,b,c:\Re^+ \times \Re^n \to \Re$ with $a(t,0)=b(t,0)=c(t,0)=0$ for all $t \ge 0$, a function $q \in \mathcal{E}$ and a $C^0$ positive definite function $\rho:\Re^+ \to \Re^+$. We showed in Example 3.3 that $V:\Re^+ \times \Re^n \to \Re^+$ is an ORCLF for (1.1) provided that implications (3.9), (3.10) hold. In this example we show that implication (3.10) only is sufficient to guarantee that $V:\Re^+ \times \Re^n \to \Re^+$ is an ORCLF for (1.1). Indeed, let $\Omega_0 := \{(t,x)\in \Re^+ \times \Re^n : a(t,x) \ge 0\}$ and $\Omega_1 := \{(t,x)\in \Re^+ \times \Re^n : a(t,x) < 0\}$. Moreover, define $k_1:\Omega_1 \to \Re$ by

$$k_1(t,x) := \frac{\sqrt{b^2(t,x) + 4|a(t,x)c(t,x)| + 4|a(t,x)|\rho(V(t,x))} - b(t,x)}{2a(t,x)}\,,\ (t,x)\in\Omega_1 \tag{3.15}$$

The specification of the set-valued map $\mathcal{U}(t,x,u)$ for $(t,x)\in\Omega_0$ has been given in Example 3.3. Therefore, the function $\Psi:\Re^+\times\Re^n\times U\to\Re\cup\{+\infty\}$ defined by:

$$\Psi(t,x,u):=\rho\big(V(t,x)\big)+a(t,x)u^2+b(t,x)u+c(t,x)\text{, for }(t,x)\in\Omega_0 \tag{3.16a}$$

$$\Psi(t,x,u):=\rho\big(V(t,x)\big)+\sup\left\{\,a(t,x)v^2+b(t,x)v+c(t,x):\ |v-k_1(t,x)|\le|u-k_1(t,x)|\,\right\}\text{, for }(t,x)\in\Omega_1 \tag{3.16b}$$

Clearly, $\Psi(t,x,u)$ as defined by (3.16a) is continuous on the interior of $\Omega_0$. Furthermore, it follows from continuity of $a(t,x), b(t,x), c(t,x), k_1(t,x)$ on $\Omega_1$ and Theorem 1.4.16 in [4] that $\Psi(t,x,u)$ as defined by (3.16b) is upper semi-continuous on $\Omega_1$. The reader should notice that implication (3.10) guarantees that $\Psi:\Re^+\times\Re^n\times U\to\Re\cup\{+\infty\}$ is upper semi-continuous (since $\Psi(t,x,u)\le\rho\big(V(t,x)\big)-\dfrac{b^2(t,x)}{4a(t,x)}+c(t,x)$ for all $(t,x,u)\in\Omega_1\times\Re$). Consequently, Lemma 3.4 guarantees that property (ii) of Definition 2.6 holds with $\Psi:\Re^+\times\Re^n\times U\to\Re\cup\{+\infty\}$ as defined by (3.16) and that $V:\Re^+\times\Re^n\to\Re^+$ is an ORCLF for (1.1). $\triangleleft$

## 4. Extensions to Systems Described by Retarded Functional Differential Equations

In this section we extend the methodology presented in Section 2, to infinite-dimensional systems described by Retarded Functional Differential Equations (RFDEs). Particularly, we consider control systems of the form (1.3) under the following hypotheses:

**(S1)** The mapping $(x,u,d)\to f(t,d,x,u)$ is continuous for each fixed $t\ge 0$ and such that for every bounded $I\subseteq\Re^+$ and for every bounded $S\subset C^0([-r,0];\Re^n)\times U$, there exists a constant $L\ge 0$ such that:

$$\big(x(0)-y(0)\big)'\big(f(t,d,x,u)-f(t,d,y,u)\big)\le L\max_{\tau\in[-r,0]}|x(\tau)-y(\tau)|^2=L\|x-y\|_r^2$$
$$\forall t\in I\,,\ \forall(x,u,y,u)\in S\times S\,,\ \forall d\in D$$

Hypothesis (S1) is equivalent to the existence of a continuous non-decreasing function $L:\Re^+\to\Re^+$, with the following property:

$$\begin{aligned}&\big(x(0)-y(0)\big)'\big(f(t,d,x,u)-f(t,d,y,u)\big)\le L(t+\|x\|_r+\|y\|_r+|u|)\,\|x-y\|_r^2\\&\forall(t,x,y,d,u)\in\Re^+\times C^0([-r,0];\Re^n)\times C^0([-r,0];\Re^n)\times D\times U\end{aligned} \tag{4.1}$$

**(S2)** For every bounded $\Omega\subset\Re^+\times D\times C^0([-r,0];\Re^n)\times U$ the image set $f(\Omega)\subset\Re^n$ is bounded.

**(S3)** There exists a countable set $A\subset\Re^+$, which is either finite or $A=\{t_k\ ;\ k=1,...,\infty\}$ with $t_{k+1}>t_k>0$ for all $k=1,2,...$ and $\lim t_k=+\infty$, such that mapping $(t,x,u,d)\in(\Re^+\setminus A)\times C^0([-r,0];\Re^n)\times U\times D\to f(t,d,x,u)$ is continuous. Moreover, for each fixed $(t_0,x,u,d)\in\Re^+\times C^0([-r,0];\Re^n)\times U\times D$, we have $\lim_{t\to t_0^+}f(t,d,x,u)=f(t_0,d,x,u)$.

**(S4)** For every $\varepsilon>0$, $t\in\Re^+$, there exists $\delta:=\delta(\varepsilon,t)>0$ such that $\sup\left\{\,|f(\tau,d,x,u)|\ ;\ \tau\in\Re^+\ ,\ d\in D, u\in U\ ,\ |\tau-t|+\|x\|_r+|u|<\delta\,\right\}<\varepsilon$.

**(S5)** The mapping $u\to f(t,d,x,u)$ is Lipschitz on bounded sets, in the sense that for every bounded $I\subseteq\Re^+$ and for every bounded $S\subset C^0([-r,0];\Re^n)\times U$, there exists a constant $L_U\ge 0$ such that:

$$|f(t,d,x,u)-f(t,d,x,v)|\le L_U|u-v|\,,\ \forall t\in I\,,\ \forall(x,u,x,v)\in S\times S\,,\ \forall d\in D$$

Hypothesis (S5) is equivalent to the existence of a continuous, non-decreasing function $L_U:\Re^+\to\Re^+$, with the following property:

$$\left| f(t,d,x,u) - f(t,d,x,v) \right| \le L_U\,(t + \|x\|_r + |u| + |v|)\,|u - v|$$
$$\forall (t,x,d,u,v) \in \Re^+ \times C^0([-r,0];\Re^n) \times D \times U \times U \tag{4.2}$$

**(S6):** The set $D \subset \Re^l$ is compact and $U \subseteq \Re^m$ is a closed convex set.

**(S7)** The mapping $H(t,x)$ is Lipschitz on bounded sets, in the sense that for every bounded $I \subseteq \Re^+$ and for every bounded $S \subset C^0([-r,0];\Re^n)$, there exists a constant $L_H \ge 0$ such that:

$$\|H(t,x) - H(\tau,y)\|_Y \le L_H\left(|t-\tau| + \|x-y\|_r\right),\ \forall (t,\tau) \in I \times I\ ,\ \forall (x,y) \in S \times S$$

Following the methodology described in Section 2, we next analyze in detail the following steps of the method:

**4.I.** Notions of Output Stability
**4.II.** Lyapunov-like criteria for Output stability
**4.III.** Definition of the Output Robust Control Lyapunov Functional
**4.IV.** Converse Lyapunov theorems for output stability

## 4.I. Notions of Output Stability

We consider uncertain dynamical systems described by RFDEs of the form:

$$\begin{aligned} &\dot{x}(t) = f(t, d(t), T_r(t)x)\ ,\ t \ge t_0 \\ &Y(t) = H(t, T_r(t)x) \\ &x(t) \in \Re^n\ , Y(t) \in Y,\ d(t) \in D \end{aligned} \tag{4.3}$$

where $r > 0$ is a constant, $f : \Re^+ \times D \times C^0([-r,0];\Re^n) \to \Re^n$, $H : \Re^+ \times C^0([-r,0];\Re^n) \to Y$ satisfy $f(t,d,0) = 0$, $H(t,0) = 0$ for all $(t,d) \in \Re^+ \times D$, $D \subseteq \Re^l$ is a non-empty compact set, $Y$ is a normed linear space and $T_r(t)x = x(t+\theta)$ ; $\theta \in [-r,0]$, under the following hypotheses:

**(Q1)** The mapping $(x,d) \to f(t,d,x)$ is continuous for each fixed $t \ge 0$ and such that for every bounded $I \subseteq \Re^+$ and for every bounded $S \subset C^0([-r,0];\Re^n)$, there exists a constant $L \ge 0$ such that

$$\begin{aligned} &\left(x(0) - y(0)\right)'\left(f(t,d,x) - f(t,d,y)\right) \le L\|x - y\|_r^2 \\ &\forall t \in I\ ,\ \forall (x,y) \in S \times S\ ,\ \forall d \in D \end{aligned}$$

**(Q2)** For every bounded $\Omega \subset \Re^+ \times D \times C^0([-r,0];\Re^n)$ the image set $f(\Omega) \subset \Re^n$ is bounded.

**(Q3)** There exists a countable set $A \subset \Re^+$, which is either finite or $A = \{t_k\ ;\ k = 1,...,\infty\}$ with $t_{k+1} > t_k > 0$ for all $k = 1,2,...$ and $\lim t_k = +\infty$, such that mapping $(t,x,d) \in (\Re^+ \setminus A) \times C^0([-r,0];\Re^n) \times D \to f(t,d,x)$ is continuous. Moreover, for each fixed $(t_0,x,d) \in \Re^+ \times C^0([-r,0];\Re^n) \times D$, we have $\lim\limits_{t \to t_0^+} f(t,d,x) = f(t_0,d,x)$.

**(Q4)** For every $\varepsilon > 0$, $t \in \Re^+$, there exists $\delta := \delta(\varepsilon,t) > 0$ such that $\sup\left\{ |f(\tau,d,x)|\ ;\ \tau \in \Re^+\ ,\ d \in D\ ,\ |\tau - t| + \|x\|_r < \delta \right\} < \varepsilon$.

**(Q5)** Hypothesis (S7) holds for the output map.

It should be emphasized for systems of the form (4.3) under hypotheses (Q1-5) that

- $0 \in C^0([-r,0];\Re^n)$ is a robust equilibrium point in the sense described in [11,15] for system (4.3) under hypotheses (Q1-5),

- system (4.3) under hypotheses (Q1-5) satisfies the "Boundedness-Implies-Continuation" property and the classical semigroup property (see [11,15]),

- for each $(t_0, x_0, d) \in \Re^+ \times C^0([-r,0];\Re^n) \times M_D$ there exists $t_{\max} \in (t_0,+\infty]$ and a unique continuous mapping $x:[t_0 - r, t_{\max}) \to \Re^n$ (the solution of (4.3)) being absolutely continuous on $[t_0, t_{\max})$ with $x(t_0+\theta) = x_0(\theta)$ for all $\theta \in [-r,0]$ and $\dot{x}(t) = f(t, d(t), T_r(t)x)$ a.e. for $t \in [t_0, t_{\max})$. Moreover, if $t_{\max} < +\infty$ then $\limsup_{t \to t_{\max}^-} |x(t)| = +\infty$. We denote the solution of (4.3) with initial condition $T_r(t_0)x = x_0$ corresponding to $d \in M_D$ by $x(t, t_0, x_0, d)$ and by $Y(t, t_0, x_0, d)$ we denote the output of system (4.3), i.e., $Y(t, t_0, x_0, d) = H(t, x(t, t_0, x_0, d))$.

For systems of the form (4.3) under hypotheses (Q1-5) we adopt the definitions of RGAOS and URGAOS given in [15] for a wide class of deterministic systems with disturbances. For completeness we repeat the definitions here.

**Definition 4.1:** *We say that (4.3) under hypotheses (Q1-5) is* ***Robustly Forward Complete (RFC)*** *if for every $s \ge 0$, $T \ge 0$, it holds that*

$$\sup\left\{ \left\| x(t_0+\xi, t_0, x_0, d) \right\|_r \; ; \; \xi \in [0,T], \left\| x_0 \right\|_r \le s, t_0 \in [0,T], d \in M_D \right\} < +\infty$$

**Definition 4.2:** *Consider system (4.3) under hypotheses (Q1-5). We say that (4.3) is* ***Robustly Globally Asymptotically Output Stable (RGAOS),*** *if (4.3) is RFC and the following properties hold:*

***P1(Output Stability)*** *For every $\varepsilon > 0$, $T \ge 0$, it holds that*

$$\sup\left\{ \left\| Y(t, t_0, x_0, d) \right\|_Y ; t \ge t_0 , \left\| x_0 \right\|_r \le \varepsilon , t_0 \in [0,T], d \in M_D \right\} < +\infty$$

*and there exists a $\delta := \delta(\varepsilon, T) > 0$ such that*

$$\left\| x_0 \right\|_r \le \delta , t_0 \in [0,T] \Rightarrow \left\| Y(t, t_0, x_0, d) \right\|_Y \le \varepsilon , \forall t \ge t_0 , \forall d \in M_D$$

***P2(Uniform Output Attractivity on bounded sets of initial data)*** *For every $\varepsilon > 0$, $T \ge 0$ and $R \ge 0$, there exists a $\tau := \tau(\varepsilon, T, R) \ge 0$, such that*

$$\left\| x_0 \right\|_r \le R , t_0 \in [0,T] \Rightarrow \left\| Y(t, t_0, x_0, d) \right\|_Y \le \varepsilon , \forall t \ge t_0 + \tau , \forall d \in M_D$$

**Definition 4.3:** *Consider system (4.3) under hypotheses (Q1-5). We say that (4.3) is Uniformly* ***Robustly Globally Asymptotically Output Stable (URGAOS),*** *if (4.3) is RFC and the following properties hold:*

***P1(Uniform Output Stability)*** *For every $\varepsilon > 0$, it holds that*

$$\sup\left\{ \left\| Y(t, t_0, x_0, d) \right\|_Y ; t \ge t_0 , \left\| x_0 \right\|_r \le \varepsilon , t_0 \ge 0 , d \in M_D \right\} < +\infty$$

*and there exists a $\delta := \delta(\varepsilon) > 0$ such that*

$$\left\| x_0 \right\|_r \le \delta , t_0 \ge 0 \Rightarrow \left\| Y(t, t_0, x_0, d) \right\|_Y \le \varepsilon , \forall t \ge t_0 , \forall d \in M_D$$

***P2(Uniform Output Attractivity on bounded sets of initial states)*** *For every $\varepsilon > 0$ and $R \ge 0$, there exists a $\tau := \tau(\varepsilon, R) \ge 0$, such that*

$$\left\| x_0 \right\|_r \le R , t_0 \ge 0 \Rightarrow \left\| Y(t, t_0, x_0, d) \right\|_Y \le \varepsilon , \forall t \ge t_0 + \tau , \forall d \in M_D$$

Obviously, the notions of RGAOS, URGAOS are direct extensions of the notions of RGAOS and URGAOS for finite-dimensional systems.

## 4.II. Lyapunov-like criteria for Output stability

Let $x \in C^0\left([-r,0];\Re^n\right)$ and $V:\Re^+ \times C^0\left([-r,0];\Re^n\right) \to \Re$ be a locally bounded functional. By $E_h(x;v)$, where $0 \le h < r$ and $v \in \Re^n$ we denote the following operator:

$$E_h(x;v) := \begin{cases} x(0)+(\theta+h)v & \text{for } \; -h<\theta\le 0 \\ x(\theta+h) & \text{for } \; -r\le\theta\le -h \end{cases} \tag{4.4}$$

and we define

$$V^0(t,x;v) := \limsup_{\substack{h\to 0^+ \\ y\to 0,\, y\in C^0([-r,0];\Re^n)}} \frac{V(t+h, E_h(x;v)+hy)-V(t,x)}{h} \tag{4.5}$$

The following lemma presents some elementary properties of the generalized derivative given above. Notice that the function $(t,x,v)\to V^0(t,x;v)$ may take values in the extended real number set $\Re^* = [-\infty,+\infty]$. Its proof is almost identical with Lemma 2.7 in [14] and is omitted.

**Lemma 4.4:** *Let $V:\Re^+ \times C^0\left([-r,0];\Re^n\right) \to \Re$ be a locally bounded functional and let $x \in C^0([t_0-r,t_{\max});\Re^n)$ a solution of (4.3) under hypotheses (Q1-5) with initial condition $x(t_0)=x_0 \in C^0([-r,0];\Re^n)$, corresponding to certain $d \in M_D$, where $t_{\max}\in(t_0,+\infty]$ is the maximal existence time of the solution. Then it holds that*

$$\limsup_{h\to 0^+} h^{-1}\left(V(t+h,T_r(t+h)x)-V(t,T_r(t)x)\right) \le V^0(t,T_r(t)x;D^+x(t)) \text{ , a.e. on } [t_0,t_{\max}) \tag{4.6}$$

*where $D^+x(t) = \lim_{h\to 0^+} h^{-1}\left(x(t+h)-x(t)\right)$.*

An important class of functionals is presented next.

**Definition 4.5:** *We say that a continuous functional $V:\Re^+\times C^0([-r,0];\Re^n)\to\Re^+$, is "almost Lipschitz on bounded sets", if there exist non-decreasing functions $L_V:\Re^+\to\Re^+$, $P:\Re^+\to\Re^+$, $G:\Re^+\to[1,+\infty)$ such that for all $R\ge 0$, the following properties hold:*

**(P1)** *For every $x,y \in \left\{x\in C^0([-r,0];\Re^n)\,;\, \|x\|_r \le R\right\}$, it holds that:*

$$\left|V(t,y)-V(t,x)\right| \le L_V(R)\,\|y-x\|_r \,, \;\; \forall t\in[0,R]$$

*(i.e., the mapping $\Re^+\times C^0([-r,0];\Re^n) \ni (t,x) \to V(t,x)\in\Re^+$ is completely locally Lipschitz with respect to $x\in C^0([-r,0];\Re^n)$)*

**(P2)** *For every absolutely continuous function $x:[-r,0]\to\Re^n$ with $\|x\|_r \le R$ and essentially bounded derivative, it holds that:*

$$\left|V(t+h,x)-V(t,x)\right| \le hP(R)\left(1+\sup_{-r\le\tau\le 0}\left|\dot{x}(\tau)\right|\right), \text{ for all } t\in[0,R] \text{ and } 0\le h\le \frac{1}{G\left(R+\sup_{-r\le\tau\le 0}\left|\dot{x}(\tau)\right|\right)}$$

The reader should notice that for functionals $V:\Re^+\times C^0([-r,0];\Re^n)\to\Re^+$, which are is almost Lipschitz on bounded sets we obtain the following simplification for the derivative $V^0(t,x;v)$ defined by (4.5) for all $(t,x,v)\in\Re^+\times C^0\left([-r,0];\Re^n\right)\times\Re^n$:

$$V^0(t,x;v) = \limsup_{h \to 0^+} \frac{V(t+h, E_h(x;v)) - V(t,x)}{h}$$

The following proposition is based on the results obtained in [16] and provides Lyapunov-like criteria for RGAOS and URGAOS for (4.3). Its proof is provided in the Appendix.

**Proposition 4.6:** *Consider system (4.3) under hypotheses (Q1-5). Suppose that there exist functions* $a_1, a_2 \in K_\infty$, $\beta, \mu \in K^+$, $q \in \mathcal{E}$, *a positive definite continuous function* $\rho : \Re^+ \to \Re^+$ *and a mapping* $V : \Re^+ \times C^0([-r,0];\Re^n) \to \Re^+$, *which is almost Lipschitz on bounded sets, such that the following inequalities hold for all* $(t,x,d) \in \Re^+ \times C^0([-r,0];\Re^n) \times D$ *:*

$$\max\left\{ a_1\left( \|H(t,x)\|_Y \right), a_1\left( \mu(t)\|x\|_r \right) \right\} \le V(t,x) \le a_2\left( \beta(t)\|x\|_r \right) \tag{4.7}$$

$$V^0(t,x; f(t,d,x)) \le -\rho\left(V(t,x)\right) + q(t) \tag{4.8}$$

*Then system (4.3) is RGAOS. Moreover, if* $\beta(t) \equiv 1$ *and* $q(t) \equiv 0$ *then system (4.3) is URGAOS.*

## 4.III. Definition of the Output Robust Control Lyapunov Functional

We next give the definition of the Output Robust Control Lyapunov Functional for system (1.3). The definition is in the same spirit with Definition 2.6 of the notion of ORCLF for finite-dimensional control systems.

**Definition 4.7:** *We say that (1.3) admits an* ***Output Robust Control Lyapunov Functional (ORCLF)*** *if there exists an almost Lipschitz on bounded sets functional* $V : \Re^+ \times C^0([-r,0];\Re^n) \to \Re^+$ *(called the Output Control Lyapunov Functional), which satisfies the following properties:*

**(i)** *There exist functions* $a_1, a_2 \in K_\infty$, $\beta, \mu \in K^+$ *such that inequality (4.7) holds for all* $(t,x) \in \Re^+ \times C^0\left([-r,0];\Re^n\right)$.

**(ii)** *There exists a function* $\Psi : \Re^+ \times \Re^q \times U \to \Re \cup \{+\infty\}$ *with* $\Psi(t,0,0) = 0$ *for all* $t \ge 0$ *such that for each* $u \in U$ *the mapping* $(t,\varphi) \to \Psi(t,\varphi,u)$ *is upper semi-continuous, a function* $q \in \mathcal{E}$, *a continuous mapping* $\Re^+ \times C^0([-r,0];\Re^n) \ni (t,x) \to \Phi(t,x) \in \Re^p$ *being completely locally Lipschitz with respect to* $x \in C^0([-r,0];\Re^n)$ *with* $\Phi(t,0) = 0$ *for all* $t \ge 0$ *and a* $C^0$ *positive definite function* $\rho : \Re^+ \to \Re^+$ *such that the following inequality holds:*

$$\inf_{u \in U} \Psi(t,\varphi,u) \le q(t)\,,\ \forall t \ge 0, \forall \varphi = (\varphi_1,...,\varphi_p)' \in \Re^p \tag{4.9}$$

*Moreover, for every finite set* $\{u_1, u_2, ..., u_N\} \subset U$ *and for every* $\lambda_i \in [0,1]$ *(* $i = 1,...,N$ *) with* $\sum_{i=1}^{N} \lambda_i = 1$, *it holds that:*

$$\sup_{d \in D} V^0\left( t,x; f\left( t,d,x,\sum_{i=1}^{N} \lambda_i\, u_i \right) \right) \le -\rho\left(V(t,x)\right) + \max\left\{ \Psi(t,\Phi(t,x),u_i)\,, i = 1,...,N \right\},$$
$$\forall (t,x) \in \Re^+ \times C^0\left([-r,0];\Re^n\right) \tag{4.10}$$

*If in addition to the above there exist* $a \in K_\infty$, $\gamma \in K^+$ *such that for every* $(t,\varphi) \in \Re^+ \times \Re^q$ *there exists* $u \in U$ *with* $|u| \le a\left(\gamma(t)|\varphi|\right)$ *such that*

$$\Psi(t,\varphi,u) \le q(t) \tag{4.11}$$

*then we say that* $V : \Re^+ \times C^0([-r,0];\Re^n) \to \Re^+$ *satisfies the "small-control" property.*

*For the case* $H(t,x) \equiv x \in C^0([-r,0];\Re^n)$ *we simply call* $V : \Re^+ \times C^0([-r,0];\Re^n) \to \Re^+$ *a State Robust Control Lyapunov Functional (SRCLF).*

**Remark 4.8:** It should be clear that in the finite-dimensional case the continuous mapping $\Phi(t,x)=(\Phi_1(t,x),...,\Phi_p(t,x))'$ is replaced by the mapping $\Phi(t,x):=x\in\Re^n$ with $p=n$. The question of the construction of the mapping $\Psi:\Re^+\times\Re^p\times U\to\Re\cup\{+\infty\}$ can be handled with exactly the same way as shown in Section 3, provided that we can find appropriate continuous mappings $\Phi(t,x)=(\Phi_1(t,x),...,\Phi_p(t,x))'$, $G:\Re^+\times\Re^p\times U\to\Re$ with $\rho\big(V(t,x)\big)+\sup_{d\in D}V^0\big(t,x;\,f\big(t,d,x,u\big)\ \big)\le G(t,\Phi(t,x),u)$ for all $(t,x)\in\Re^+\times C^0([-r,0];\Re^n)$ and $\inf_{u\in U}G(t,\varphi,u)\le q(t)$ for all $(t,\varphi)\in\Re^+\times\Re^p$. In this case all constructions of $\Psi:\Re^+\times\Re^p\times U\to\Re\cup\{+\infty\}$ given in Section 3 may be repeated with the quantity $\frac{\partial V}{\partial t}(t,x)+\sup_{d\in D}\left(\frac{\partial V}{\partial x}(t,x)f(t,d,x,u)\ \right)+\rho\big(V(t,x)\big)$ replaced by the quantity $G(t,\varphi,u)$.

## 4.IV. Converse Lyapunov theorems for output stability

In this work we are going to exploit the converse Lyapunov theorems for RGAOS and URGAOS presented in [16].

## 4.V. Main Results

We are now in a position to state and prove our main results for the infinite-dimensional case (1.3).

**Theorem 4.9:** *Consider system (1.3) under hypotheses (S1-7). The following statements are equivalent:*

**(a)** *There exists a continuous mapping $\Re^+\times C^0([-r,0];\Re^n)\ni(t,x)\to k(t,x)\in U$ being completely locally Lipschitz with respect to $x\in C^0([-r,0];\Re^n)$ with $k(t,0)=0$ for all $t\ge 0$, such that the closed-loop system (1.3) with $u=k(t,T_r(t)x)$ is RGAOS.*

**(b)** *System (1.3) admits an ORCLF, which satisfies the small control property with $q(t)\equiv 0$.*

**(c)** *System (1.3) admits an ORCLF.*

**Theorem 4.10:** *Consider system (1.3) under hypotheses (S1-7). The following statements are equivalent:*

**(a)** *System (1.3) admits an ORCLF, which satisfies the small-control property and inequalities (4.7), (4.11) with $\beta(t)\equiv 1$, $q(t)\equiv 0$. Moreover, there exist continuous mappings $\eta\in K^+$, $A\ni(t,\varphi)\to K(t,\varphi)\in U$ where $A=\bigcup_{t\ge 0}\{t\}\times\{\varphi\in\Re^p:|\varphi|<4\eta(t)\}$ being locally Lipschitz with respect to $\varphi$ with $K(t,0)=0$ for all $t\ge 0$ and such that*

$$\Psi(t,\Phi(t,x),K(t,\Phi(t,x)))\le 0,\ \text{for all } (t,x)\in\Re^+\times C^0([-r,0];\Re^n) \text{ with } |\Phi(t,x)|\le 2\eta(t) \tag{4.12}$$

*where $\Phi=(\Phi_1,...,\Phi_p)':\Re^+\times C^0([-r,0];\Re^n)\to\Re^p$ and $\Psi:\Re^+\times\Re^p\times U\to\Re\cup\{+\infty\}$ are the mappings involved in property (ii) of Definition 4.7.*

**(b)** *There exists a continuous mapping $\Re^+\times C^0([-r,0];\Re^n)\ni(t,x)\to k(t,x)\in U$ being completely locally Lipschitz with respect to $x\in C^0([-r,0];\Re^n)$ with $k(t,0)=0$ for all $t\ge 0$, such that the closed-loop system (1.3) with $u=k(t,T_r(t)x)$ is URGAOS.*

**Remark 4.11:** From the proof of Theorem 4.10 it will become apparent that if statement (a) of Theorem 4.10 is strengthened so that the ORCLF $V$, the mappings $\Phi=(\Phi_1,...,\Phi_p)':\Re^+\times C^0([-r,0];\Re^n)\to\Re^p$, $\Psi$ involved in property (ii) of Definition 4.7 and the mapping $K:A\to U$ are time independent then the continuous mapping $k$, whose existence is guaranteed by statement (b) of Theorem 4.10, is time invariant.

**Proof of Theorem 4.9:** The implication (b) $\Rightarrow$ (c) is obvious and we prove implications (a) $\Rightarrow$ (b) and (c) $\Rightarrow$ (a).

**(c) $\Rightarrow$ (a)** Suppose that (1.3) admits an ORCLF. Without loss of generality, we may assume that the function $q \in \mathcal{E}$ involved in (4.9) is positive for all $t \geq 0$.

Furthermore define:

$$\Xi(t,\varphi,u) := \Psi(t,\varphi,u) - 8q(t)\,,\ (t,\varphi,u) \in \Re^+ \times \Re^p \times U \tag{4.13a}$$

$$\Xi(t,\varphi,u) := \Xi(0,\varphi,u)\,,\ (t,\varphi,u) \in (-1,0) \times \Re^p \times U \tag{4.13b}$$

The definition of $\Xi$, given by (4.13a,b), guarantees that the function $\Xi : (-1,+\infty) \times \Re^p \times U \to \Re \cup \{+\infty\}$ with $\Xi(t,0,0) = -8q(\max\{0,t\})$ for all $t > -1$ is such that, for each $u \in U$ the mapping $(t,\varphi) \to \Xi(t,\varphi,u)$ is upper semi-continuous. By virtue of (4.9) and upper semi-continuity of $\Xi$, it follows that for each $(t,\varphi) \in (-1,+\infty) \times \Re^p$ there exist $u = u(t,\varphi) \in U$ and $\delta = \delta(t,\varphi) \in (0, t+1)$, such that

$$\Xi(\tau, y, u(t,\varphi)) \leq 0\,,\ \forall (\tau, y) \in \left\{(\tau, y) \in (-1,+\infty) \times \Re^p : |\tau - t| + |y - \varphi| < \delta\right\} \tag{4.14}$$

Using (4.14) and standard partition of unity arguments, we can determine sequences $\{(t_i, \varphi_i) \in (-1,+\infty) \times \Re^p\}_{i=1}^{\infty}$, $\{u_i \in U\}_{i=1}^{\infty}$, $\{\delta_i\}_{i=1}^{\infty}$ with $\delta_i = \delta(t_i, \varphi_i) \in (0, t_i + 1)$ associated with a sequence of open sets $\{\Omega_i\}_{i=1}^{\infty}$ with

$$\Omega_i \subseteq \left\{(\tau, y) \in (-1,+\infty) \times \Re^p : |\tau - t_i| + |y - \varphi_i| < \delta_i\right\} \tag{4.15}$$

forming a locally finite open covering of $(-1,+\infty) \times \Re^p$ and in such a way that:

$$\Xi(\tau, y, u_i) \leq 0\,,\ \forall (\tau, y) \in \Omega_i \tag{4.16}$$

Also, a family of smooth functions $\{\theta_i\}_{i=1}^{\infty}$ with $\theta_i(t,\varphi) \geq 0$ for all $(t,x) \in (-1,+\infty) \times \Re^p$ can be determined with

$$supp\ \theta_i \subseteq \Omega_i \tag{4.17}$$

$$\sum_{i=1}^{\infty} \theta_i(t,\varphi) = 1\,,\ \forall (t,\varphi) \in (-1,+\infty) \times \Re^p \tag{4.18}$$

Using exactly the same methodology as in the proof of implication (d) $\Rightarrow$ (a) of Theorem 2.8 and the facts that $\Xi(t,0,0) = -8q(t) < 0$ for all $t \geq 0$ and that the mapping $(t,\varphi) \to \Xi(t,\varphi,0)$ is upper semi-continuous, we may establish the existence of a $C^\infty$ positive function $\eta : \Re^+ \to (0,+\infty)$ with the following property:

$$|\varphi| \leq 2\eta(t) \ \Rightarrow \Xi(t,\varphi,0) \leq 0 \tag{4.19}$$

Let $h \in C^\infty(\Re;[0,1])$ be a smooth non-decreasing function with $h(s) = 0$ for all $s \leq 0$ and $h(s) = 1$ for all $s \geq 1$. We define for all $(t,x) \in \Re^+ \times C^0\left([-r,0];\Re^n\right)$:

$$k(t,x) := h\left(\frac{|\Phi(t,x)|^2 - 2\eta^2(t)}{2\eta^2(t)}\right) \sum_{i=1}^{\infty} \theta_i(t, \Phi(t,x)) u_i \tag{4.20}$$

where $\Phi = (\Phi_1, ..., \Phi_p)' : \Re^+ \times C^0([-r,0];\Re^n) \to \Re^p$ is the mapping involved in property (ii) of Definition 4.7. Clearly, $k$ as defined by (4.20) is a mapping satisfying the property that for every bounded $\Omega \subset \Re^{+} \times C^{0}([-r,0];\Re^{n})$ it holds that $\sup\left\{ \frac{|k(t,x)-k(t,y)|}{\|x-y\|_{r}} : (t,x)\in\Omega, (t,y)\in\Omega, x\neq y \right\} < +\infty$, with $k(t,0)=0$ for all $t \geq 0$. Moreover, since $k(t,x)$ is defined as a (finite) convex combination of $u_i \in U$ and $0 \in U$, we have $k(t,x)\in U$ for all $(t,x)\in \Re^{+} \times C^{0}\left([-r,0];\Re^{n}\right)$.

Let $(t,x)\in \Re^{+} \times C^{0}\left([-r,0];\Re^{n}\right)$ with $|\Phi(t,x)| \geq 2\,\eta(t)$ and define $J(t,x)=\left\{ j\in\{1,2,...\}\,;\, \theta_j(t,\Phi(t,x)) \neq 0 \right\}$ (a finite set). Notice that by virtue of (4.10) and definition (4.20) we get:

$$\begin{aligned} &\sup_{d\in D} V^{0}\left( t,x;\, f\left(t,d,x,k(t,x)\right)\right) = \sup_{d\in D} V^{0}\left( t,x;\, f\left( t,d,x, \sum_{j\in J(t,x)} \theta_j(t,\Phi(t,x))u_j \right)\right) \\ &\leq -\rho\left(V(t,x)\right) + \max_{j\in J(t,x)} \left\{ \Psi(t,\Phi(t,x),u_j) \right\} \end{aligned} \tag{4.21}$$

Notice that for each $j\in J(t,x)$ we obtain from (4.17) that $(t,\Phi(t,x))\in\Omega_j$. Consequently, by virtue of (4.16) and definition (4.13a) we have that $\Psi(t,\Phi(t,x),u_j) \leq 8q(t)$, for all $j\in J(t,x)$. Combining the previous inequality with inequality (4.21), we conclude that the following property holds for all $(t,x,d)\in \Re^{+} \times C^{0}\left([-r,0];\Re^{n}\right)\times D$ with $|\Phi(t,x)| \geq 2\,\eta(t)$:

$$V^{0}(t,x;\, f(t,d,x,k(t,x))) \leq -\rho\left(V(t,x)\right) + 8q(t) \tag{4.22}$$

Let $(t,x)\in \Re^{+} \times C^{0}\left([-r,0];\Re^{n}\right)$ with $|\Phi(t,x)| \leq \sqrt{2}\,\eta(t)$. Notice that by virtue of definition (4.20) we get:

$$\sup_{d\in D} V^{0}\left( t,x;\, f\left(t,d,x,k(t,x)\right)\right) = \sup_{d\in D} V^{0}\left( t,x;\, f\left(t,d,x,0\right)\right)$$

By virtue of (4.10), (4.19), (4.13a) and the above inequality we conclude that (4.22) holds as well for all $(t,x)\in \Re^{+} \times C^{0}\left([-r,0];\Re^{n}\right)$ with $|\Phi(t,x)| \leq \sqrt{2}\,\eta(t)$. Finally, for the case $(t,x)\in \Re^{+} \times C^{0}\left([-r,0];\Re^{n}\right)$ with $\sqrt{2}\;\eta(t) < |\Phi(t,x)| < 2\eta(t)$, let $J(t,x)=\left\{ j\in\{1,2,...\}\,;\, \theta_j(t,\Phi(t,x)) \neq 0 \right\}$ and notice that from (4.10) we get:

$$\begin{aligned} &\sup_{d\in D} V^{0}\left( t,x;\, f\left(t,d,x,k(t,x)\right)\right) = \sup_{d\in D} V^{0}\left( t,x;\, f\left( t,d,x, h\left( \frac{|x|^{2} - 2\eta^{2}(t)}{2\eta^{2}(t)} \right) \sum_{j\in J(t,x)} \theta_j(t,x)u_j \right)\right) \\ &\leq -\rho\left(V(t,x)\right) + \max\left\{ \Psi(t,\Phi(t,x),0),\, \Psi(t,\Phi(t,x),u_j),\; j\in J(t,x) \right\} \end{aligned} \tag{4.23}$$

Taking into account definition (4.13a) and (4.16), (4.17), (4.19), (4.23), we may conclude that (4.22) holds as well for all $(t,x)\in \Re^{+} \times C^{0}\left([-r,0];\Re^{n}\right)$ with $\sqrt{2}\;\eta(t) < |\Phi(t,x)| < 2\eta(t)$. Consequently, (4.22) holds for all $(t,x)\in \Re^{+} \times C^{0}\left([-r,0];\Re^{n}\right)$.

It follows from (4.22) and Proposition 4.6 that system (1.3) with $u = k(t,T_r(t)x)$ is RGAOS.

**(a)** $\Rightarrow$ **(b)** Since system (1.3) with $u = k(t,T_r(t)x)$ is RGAOS, and since for every bounded $\Omega \subset \Re^{+} \times C^{0}([-r,0];\Re^{n})$ it holds that $\sup\left\{ \frac{|k(t,x)-k(t,y)|}{\|x-y\|_{r}} : (t,x)\in\Omega, (t,y)\in\Omega, x\neq y \right\} < +\infty$, it follows that the closed-loop system (1.3) with $u = k(t,T_r(t)x)$ satisfies hypotheses (Q1-5). Moreover, since system (1.3) with $u = k(t,T_r(t)x)$ is RGAOS, it follows from Theorem 3.6 in [11] that there exists $\mu \in K^{+}$ such that the following system

$$\begin{aligned} &\dot{x}(t) = f(t,d(t),T_r(t)x,k(t,T_r(t)x)) \\ &Y(t) = \|H(t,x)\|_{Y} + \mu(t)\|x\|_{r} \\ &x(t)\in\Re^{n},\, Y(t)\in\Re,\, d(t)\in D \end{aligned} \tag{4.24}$$

satisfies hypotheses (Q1-5) and is RGAOS. Notice that system (4.24) is the closed-loop system (1.3) with $u = k(t, T_r(t)x)$ and output defined by $Y(t) = \|H(t,x)\|_Y + \mu(t)\|x\|_r$. It follows from Theorem 3.5 in [16] that there exist functions $a_1, a_2 \in K_\infty$, $\beta \in K^+$ and a mapping $V : \Re^+ \times C^0([-r,0];\Re^n) \to \Re^+$, which is almost Lipschitz on bounded sets, such that:

$$a_1\left(\|H(t,x)\|_Y + \mu(t)\|x\|_r\right) \le V(t,x) \le a_2\left(\beta(t)\|x\|_r\right),\ \forall (t,x) \in \Re^+ \times C^0([-r,0];\Re^n) \tag{4.25}$$

$$V^0(t,x; f(t,d,x,k(t,x))) \le -V(t,x),\ \forall (t,x,d) \in \Re^+ \times C^0([-r,0];\Re^n) \times D \tag{4.26}$$

We next prove that $V$ is an ORCLF for (1.3). Obviously property (i) of Definition 4.7 is a consequence of inequality (4.25). Define for all $(t,\varphi) = (t,\varphi_1,\varphi_2)' \in \Re^+ \times \Re \times \Re^m$:

$$\Psi(t,\varphi_1,\varphi_2,u) := L_V(t + |\varphi_1|)\, L_U(t + |\varphi_1| + 2|\varphi_2| + |u - \Pr_U(\varphi_2)|)\, |u - \Pr_U(\varphi_2)| \tag{4.27}$$

and for all $(t,x) \in \Re^+ \times C^0([-r,0];\Re^n)$:

$$\Phi(t,x) := \begin{bmatrix} \|x\|_r \\ k(t,x) \end{bmatrix} \in \Re^{m+1} \tag{4.28}$$

where $L_U : \Re^+ \to \Re^+$ is the non-decreasing continuous function involved in (4.2) and $L_V : \Re^+ \to \Re^+$ is the non-decreasing function involved in property (P1) of Definition 4.5. The reader should notice that $\Phi(t,0) = 0$ for all $t \ge 0$ and that for every bounded $\Omega \subset \Re^+ \times C^0([-r,0];\Re^n)$ it holds that $\sup\left\{ \dfrac{|\Phi(t,x) - \Phi(t,y)|}{\|x-y\|_r} : (t,x) \in \Omega, (t,y) \in \Omega, x \ne y \right\} < +\infty$. Without loss of generality we may assume that $L_V : \Re^+ \to \Re^+$ is continuous as well. Convexity of the set $U \subseteq \Re^m$ implies that the mapping $\Re^m \ni \varphi_2 \to \Pr_U(\varphi_2)$ is continuous and consequently that the mapping $\Re^+ \times \Re^{m+1} \ni (t,\varphi_1,\varphi_2) \to \Psi(t,\varphi_1,\varphi_2,u)$ is continuous for each fixed $u \in U$. Notice that for every $\varphi = (t,\varphi_1,\varphi_2)' \in \Re^+ \times \Re \times \Re^m$ and every finite set $\{u_1,u_2,...,u_N\} \subset U$, $\lambda_i \in [0,1]$ ($i = 1,...,N$) with $\sum_{i=1}^N \lambda_i = 1$, definition (4.27) in conjunction with the fact $\left|\sum_{i=1}^N \lambda_i u_i - \Pr_U(\varphi_2)\right| \le \sum_{i=1}^N \lambda_i |u_i - \Pr_U(\varphi_2)| \le \max_{i=1,...,N} |u_i - \Pr_U(\varphi_2)|$ implies that:

$$\Psi(t,\varphi_1,\varphi_2,\sum_{i=1}^N \lambda_i u_i) \le \max_{i=1,...,N} \Psi(t,\varphi_1,\varphi_2,u_i) \tag{4.29}$$

By virtue of definitions (4.4), (4.5) and property (P1) of Definition 4.5 we get for all $(t,x,v,w) \in \Re^+ \times C^0\left([-r,0];\Re^n\right) \times \Re^n \times \Re^n$:

$$V^0(t,x;v) \le L_V(t + \|x\|_r)|v - w| + V^0(t,x;w) \tag{4.30}$$

Combining inequalities (4.2), (4.26) and (4.30) we obtain for all $(t,x,u) \in \Re^+ \times C^0\left([-r,0];\Re^n\right) \times U$:

$$\begin{aligned}
&\sup_{d\in D} V^0\left(t,x; f(t,d,x,u)\right) \le \sup_{d\in D} V^0\left(t,x; f(t,d,x,k(t,x))\right) + L_V(t + \|x\|_r) \sup_{d\in D} |f(t,d,x,u) - f(t,d,x,k(t,x))| \\
&\le -V(t,x) + L_V(t + \|x\|_r) L_U(t + \|x\|_r + |u| + |k(t,x)|)|u - k(t,x)| \\
&-V(t,x) + L_V(t + \|x\|_r) L_U(t + \|x\|_r + |u - k(t,x)| + 2|k(t,x)|)|u - k(t,x)|
\end{aligned}$$

The above inequality in conjunction with (4.29) and definitions (4.27), (4.28) implies that inequality (4.10) with $\rho(s) := s$ holds. Moreover, by virtue of definition (4.27), for every $(t,\varphi_1,\varphi_2) \in \Re^+ \times \Re \times \Re^m$ there exists $u \in U$ (namely $u = \Pr_U(\varphi_2)$) such that (4.11) holds with $q(t) \equiv 0$. Notice that $\left|\Pr_U(\varphi_2)\right| \le 2|\varphi_2| \le 2|\varphi|$, where $\varphi = (\varphi_1, \varphi_2) \in \Re^{m+1}$ and therefore the small control property holds with $a(s) := 2s$ and $\gamma(t) \equiv 1$.

The proof is complete. $\triangleleft$

**Proof of Theorem 4.10: (a)$\Rightarrow$ (b)** Suppose that (1.3) admits an ORCLF which satisfies the small control property with $q(t) \equiv 0$. Without loss of generality we may assume that the mapping $\Phi = (\Phi_1, ..., \Phi_p)' : \Re^+ \times C^0([-r,0];\Re^n) \to \Re^p$ involved in property (ii) of Definition 4.7 satisfies $\Phi_1(t,x) := V(t,x)$. Define:

$$\Xi(t,\varphi,u) := \Psi(t,\varphi,u) - \frac{1}{2}\rho\left(c'\varphi\right),\ (t,\varphi,u) \in \Re^+ \times \Re^p \times U \tag{4.31a}$$

$$\Xi(t,\varphi,u) := \Xi(0,\varphi,u),\ (t,\varphi,u) \in (-1,0) \times \Re^p \times U \tag{4.31b}$$

$$c = (1,0,...,0)' \in \Re^p \tag{4.31c}$$

The definition of $\Xi$, given by (4.31), guarantees that the function $\Xi : (-1,+\infty) \times \Re^p \times U \to \Re \cup \{+\infty\}$ with $\Xi(t,0,0) = 0$ for all $t > -1$ is such that, for each $u \in U$ the mapping $(t,\varphi) \to \Xi(t,\varphi,u)$ is upper semi-continuous. Let $\Theta := (-1,+\infty) \times \left\{\varphi \in \Re^p : c'\varphi \neq 0\right\}$, which is an open set. By virtue of (4.11) with $q(t) \equiv 0$ and upper semi-continuity of $\Xi$, it follows that for each $(t,\varphi) \in \Theta$ there exist $u = u(t,\varphi) \in U$ with $|u| \le a\left(\gamma(\max(0,t))|\varphi|\right)$ and $\delta = \delta(t,\varphi) \in (0, \min(1, t+1))$, $\delta(t,\varphi) \le \frac{|c'\varphi|}{2}$, such that

$$\Xi(\tau, y, u(t,\varphi)) \le 0,\ \forall (\tau, y) \in \left\{(\tau,y) \in \Theta : |\tau - t| + |y - \varphi| < \delta\right\} \tag{4.32}$$

Using (4.32) and standard partition of unity arguments, we can determine sequences $\{(t_i,\varphi_i) \in \Theta\}_{i=1}^{\infty}$, $\{u_i \in U\}_{i=1}^{\infty}$, $\{\delta_i\}_{i=1}^{\infty}$ with $|u_i| \le a\left(\gamma(\max(0,t_i))|\varphi_i|\right)$, $\delta_i = \delta(t_i,\varphi_i) \in (0, \min(1, t_i+1))$, $\delta_i = \delta(t_i,\varphi_i) \le \frac{|c'\varphi_i|}{2}$ associated with a sequence of open sets $\{\Omega_i\}_{i=1}^{\infty}$ with

$$\Omega_i \subseteq \left\{(\tau,y) \in \Theta : |\tau - t_i| + |y - \varphi_i| < \delta_i\right\} \tag{4.33}$$

forming a locally finite open covering of $\Theta$ and in such a way that:

$$\Xi(\tau, y, u_i) \le 0,\ \forall (\tau,y) \in \Omega_i \tag{4.34}$$

Also, a family of smooth functions $\{\theta_i\}_{i=1}^{\infty}$ with $\theta_i(t,\varphi) \ge 0$ for all $(t,\varphi) \in \Theta$ can be determined with

$$supp\ \theta_i \subseteq \Omega_i \tag{4.35}$$

$$\sum_{i=1}^{\infty} \theta_i(t,\varphi) = 1,\ \forall (t,\varphi) \in \Theta \tag{4.36}$$

We define for all $(t,x) \in \Re^+ \times C^0\left([-r,0];\Re^n\right)$:

$$k(t,x) := \left(1 - h\left(\frac{|\Phi(t,x)|^2 - 2\eta^2(t)}{2\eta^2(t)}\right)\right) K\left(t, \Pr_{Q(t)}(\Phi(t,x))\right) + h\left(\frac{|\Phi(t,x)|^2 - 2\eta^2(t)}{2\eta^2(t)}\right)\tilde{k}(t,x) \tag{4.37}$$

where

$$\tilde{k}(t,x) := \sum_{i=1}^{\infty} \theta_i(t,\Phi(t,x))u_i \text{ for } t \ge 0 \text{ , } x \ne 0 \tag{4.38a}$$

$$\tilde{k}(t,0) := 0 \text{ for } t \ge 0 \tag{4.38b}$$

where $h \in C^{\infty}(\Re;[0,1])$ be a smooth non-decreasing function with $h(s) = 0$ for all $s \le 0$ and $h(s) = 1$ for all $s \ge 1$ and $Q(t) = \{\varphi \in \Re^p : |\varphi| \le 3\eta(t)\}$. It follows from definition (4.37), (4.38) and that facts that the continuous mapping $\Phi$ is completely locally Lipschitz with respect to $x \in C^0([-r,0];\Re^n)$ and that the continuous mapping $\varphi \to K(t,\varphi)$ is locally Lipschitz that $k$ is completely locally Lipschitz with respect to $x \in C^0([-r,0];\Re^n)$ with $k(t,0) = 0$ for all $t \ge 0$. Moreover, it should be noticed that if the ORCLF $V$ and the function $\Psi$ involved in property (ii) of Definition 4.7 are time independent then the partition of unity arguments used above may be repeated on $\Theta := \{\varphi \in \Re^p : c'\varphi \ne 0\}$ instead of $\Theta := (-1,+\infty) \times \{\varphi \in \Re^p : c'\varphi \ne 0\}$. This implies that the constructed feedback is time invariant, provided that the mappings $\Phi = (\Phi_1,...,\Phi_p)' : \Re^+ \times C^0([-r,0];\Re^n) \to \Re^p$ and $K : A \to U$ are time independent too.

Exploiting the properties of the mappings $\Xi : (-1,+\infty) \times \Re^p \times U \to \Re \cup \{+\infty\}$, $\Psi : \Re^+ \times \Re^p \times U \to \Re \cup \{+\infty\}$, inequalities (4.12), (4.34), definitions (4.31a), (4.37) and the fact that the mapping $\Phi = (\Phi_1,...,\Phi_p)' : \Re^+ \times C^0([-r,0];\Re^n) \to \Re^p$ involved in property (ii) of Definition 4.7 satisfies $\Phi_1(t,x) := V(t,x)$, we may establish (exactly in the same way as in the proof of Theorem 4.9) the following inequality:

$$V^0(t,x; f(t,d,x,k(t,x))) \le -\frac{1}{2}\rho(V(t,x)), \ \forall (t,x,d) \in \Re^+ \times \Re^n \times D \tag{4.39}$$

The fact that system (1.3) with $u = k(t,T_r(t)x)$ is URGAOS follows directly from Proposition 4.6 and inequality (4.39).

**(b) $\Rightarrow$ (a)** Since for every bounded $\Omega \subset \Re^+ \times C^0([-r,0];\Re^n)$ it holds that $\sup\left\{\frac{|k(t,x)-k(t,y)|}{\|x-y\|_r} : (t,x) \in \Omega, (t,y) \in \Omega, x \ne y\right\} < +\infty$, it follows that the closed-loop system (1.3) with $u = k(t,T_r(t)x)$ satisfies hypotheses (Q1-5). For each $(t_0,x_0,d) \in \Re^+ \times C^0([-r,0];\Re^n) \times M_D$, we denote the solution of (1.3) with $u = k(t,T_r(t)x)$, initial condition $T_r(t_0)x = x_0$ corresponding to $d \in M_D$ by $x(t,t_0,x_0,d)$.

Since system (1.3) with $u = k(t,T_r(t)x)$ is RFC, it follows from Lemma 3.5 in [11] that there exists $\tilde{\mu} \in K^+$, $a \in K_\infty$ such that the following inequality holds for all $(t_0,x_0,d) \in \Re^+ \times C^0([-r,0];\Re^n) \times M_D$:

$$\tilde{\mu}(t)\|x\|_r \le a(\|x_0\|_r), \ \forall t \ge t_0 \tag{4.40}$$

Since system (1.3) with $u = k(t,T_r(t)x)$ is URGAOS, it follows from (4.40) that the following system

$$\begin{aligned}
&\dot{x}(t) = f(t,d(t),T_r(t)x,k(t,T_r(t)x)) \\
&Y(t) = \|H(t,x)\|_Y + \mu(t)\|x\|_r \\
&x(t) \in \Re^n \ ,Y(t) \in \Re, \ d(t) \in D
\end{aligned} \tag{4.41}$$

where $\mu(t) := \exp(-t)\tilde{\mu}(t)$ satisfies hypotheses (Q1-5) and is URGAOS. Notice that system (4.41) is the closed-loop system (1.3) with $u = k(t,T_r(t)x)$ and output defined by $Y(t) = \|H(t,x)\|_Y + \mu(t)\|x\|_r$. It follows from Theorem 3.6 in [16] that there exist functions $a_1, a_2 \in K_\infty$ and a mapping $V : \Re^+ \times C^0([-r,0];\Re^n) \to \Re^+$, which is almost Lipschitz on bounded sets, such that:

$$a_1\left(\|H(t,x)\|_Y+\mu(t)\|x\|_r\right)\le V(t,x)\le a_2\left(\|x\|_r\right),\ \forall (t,x)\in\Re^+\times C^0([-r,0];\Re^n) \tag{4.42}$$

$$V^0(t,x;\ f(t,d,x,k(t,x)))\le -V(t,x),\ \forall (t,x,d)\in\Re^+\times C^0([-r,0];\Re^n)\times D \tag{4.43}$$

The rest of proof is exactly the same with the proof of implication (a)$\Rightarrow$(b) of Theorem 4.9. The only additional thing that should be noticed is that the mappings $\Phi=(\Phi_1,...,\Phi_p)':\Re^+\times C^0([-r,0];\Re^n)\to\Re^p$, $\eta\in K^+$ and $A\ni(t,\varphi)\to K(t,\varphi)\in U$, may be selected so that (4.28) holds, $\eta(t)\equiv 1$ and $K(t,\varphi):=\Pr_U(\varphi_2)$ for all $(t,\varphi)=(t,\varphi_1,\varphi_2)'\in\Re^+\times\Re\times\Re^m$.

The proof is complete. $\triangleleft$

## 5. Applications to Triangular Time-Delay Control Systems

Our main result concerning triangular time-delay control systems of the form (1.4) is stated next. It must be compared to Theorem 5.1 in [4], which deals with the triangular finite-dimensional case.

**Theorem 5.1:** *Consider system (1.4), where* $r>0$, $D\subset\Re^l$ *is a compact set, the mappings* $f_i:\Re^+\times D\times C^0([-r,0];\Re^i)\to\Re$, $g_i:\Re^+\times D\times C^0([-r,0];\Re^i)\to\Re$ *(*$i=1,...,n$*) are continuous with* $f_i(t,d,0)=0$ *for all* $(t,d)\in\Re^+\times D$ *and each* $g_i:\Re^+\times D\times C^0([-r,0];\Re^i)\to\Re$ *(*$i=1,...,n$*) is completely locally Lipschitz with respect to* $x\in C^0([-r,0];\Re^i)$. *Suppose that there exists a function* $\varphi\in C^\infty(\Re^+;(0,+\infty))$ *being non-decreasing, such that for every* $i=1,...,n$, *it holds that:*

$$\frac{1}{\varphi\left(\|x\|_r\right)}\le g_i(t,d,x)\le\varphi\left(\|x\|_r\right),\ \forall (t,x,d)\in\Re^+\times C^0([-r,0];\Re^i)\times D \tag{5.1}$$

*Moreover, suppose that for every* $i=1,...,n$, *it holds that*

$$\sup\left\{\frac{|f_i(t,d,x)-f_i(t,d,y)|}{\|x-y\|_r}:\ (t,d)\in\Re^+\times D,\, x\in S,\, y\in S,\, x\ne y\right\}<+\infty\ ,$$
$$\text{for every bounded } S\subset C^0([-r,0];\Re^i) \tag{5.2}$$

*Then for every* $\sigma>0$ *there exist functions* $\mu_i\in C^\infty(\Re^i;(0,+\infty))$, $k_i\in C^\infty(\Re^i;\Re)$ *(*$i=1,...,n$*) with*

$$k_1(\xi_1):=-\mu_1(\xi_1)\,\xi_1 \tag{5.3a}$$

$$k_j(\xi_1,...,\xi_j):=-\mu_j(\xi_1,...,\xi_j)\left(\xi_j-k_{j-1}(\xi_1,...,\xi_{j-1})\right),\ j=2,...,n \tag{5.3b}$$

*such that the following functional:*

$$V(x):=\max_{\theta\in[-r,0]}\exp(2\sigma\theta)\left(x_1^2(\theta)+\sum_{j=2}^{n}\left|x_j(\theta)-k_{j-1}(x_1(\theta),...,x_{j-1}(\theta))\right|^2\right) \tag{5.4}$$

*is a State Robust Control Lyapunov Functional (SRCLF) for (1.4), which satisfies the "small-control" property. Moreover, the closed-loop system (1.4) with* $u(t)=k_n(x(t))$ *is URGAS. More specifically, the inequality* $V^0(x;v)\le -2\sigma V(x)$ *holds for all* $(t,x,d)\in\Re^+\times C^0([-r,0];\Re^n)\times D$ *with* $v=(f_1(t,d,x_1)+g_i(t,d,x_1)x_2(0),...,\ f_n(t,d,x)+g_n(t,d,x)k_n(x(0)))'\in\Re^n$.

**Remark 5.2:** The reader should notice that the feedback law $u(t)=k_n(x(t))$ is delay-independent. The proof of Theorem 5.1 will show that the functions $\mu_i \in C^\infty(\Re^i;(0,+\infty))$ ($i=1,...,n$) are obtained by a procedure similar to the backstepping procedure used for finite-dimensional triangular control systems. Consequently, as in the finite-dimensional case, the feedback design and the construction of the State Robust Control Lyapunov Functional proceed in parallel.

The proof of Theorem 5.1 is based on the following lemma. Its proof is provided at the Appendix. The reader should notice that Lemma 5.3 in conjunction with definition (5.3) of the SRCLF for system (1.4) indicate one of the complications mentioned in the Introduction encountered in the study of infinite-dimensional systems: although the differential equations (1.4) are affine in the control input $u \in \Re$, the derivative $V^0(x;v)$, where $v=(f_1(t,d,x_1)+g_i(t,d,x_1)x_2(0),..., f_n(t,d,x)+g_n(t,d,x)u)' \in \Re^n$, $x=(x_1,...,x_n)\in C^0([-r,0];\Re^n)$ is not affine in the control input $u \in \Re$.

**Lemma 5.3:** *Let $Q \in C^1(\Re^n;\Re^+)$, $\sigma>0$ and consider the functional $V: C^0([-r,0];\Re^n) \to \Re^+$ defined by:*

$$V(x) := \max_{\theta\in[-r,0]} \exp(2\sigma\theta) Q(x(\theta)) \tag{5.5}$$

*The functional $V: C^0([-r,0];\Re^n) \to \Re^+$ defined by (5.1), is Lipschitz on bounded sets of $C^0([-r,0];\Re^n)$ and satisfies:*

$$V^0(x;v) \le -2\sigma V(x) \text{ for all } (x,v)\in C^0([-r,0];\Re^n)\times\Re^n \text{ with } Q(x(0))<V(x) \tag{5.6a}$$

$$V^0(x;v) \le \max\left\{-2\sigma V(x), \nabla Q(x(0))v\right\} \text{ for all } (x,v)\in C^0([-r,0];\Re^n)\times\Re^n \text{ with } Q(x(0))=V(x) \tag{5.6b}$$

We are now in a position to provide the proof of Theorem 5.1.

**Proof of Theorem 5.1:** Inequality (5.2) in conjunction with the fact that $f_i(t,d,0)=0$ for all $(t,d)\in\Re^+\times D$ ($i=1,...,n$) implies the existence of a non-decreasing function $L\in C^\infty(\Re^+;(0,+\infty))$ such that for every $i=1,...,n$, it holds:

$$|f_i(t,d,x)| \le L\left(\|x\|_r\right)\|x\|_r, \quad \forall (t,x,d)\in\Re^+\times C^0([-r,0];\Re^i)\times D \tag{5.7}$$

Let $\sigma>0$ be a given number. We next define the functions $\mu_i\in C^\infty(\Re^i;(0,+\infty))$, $\gamma_i\in C^\infty(\Re^+;(0,+\infty))$, $b_i\in C^\infty(\Re^+;(0,+\infty))$ ($i=1,...,n$) using the following algorithm:

Step $i=1$: We define:

$$\mu_1(\xi_1) := \frac{\gamma_1(1+\xi_1^2)+n\sigma}{b_1(1+\xi_1^2)} \tag{5.8}$$

where $b_1\in C^\infty(\Re;(0,+\infty))$ is the function involved in (5.1) for $i=1$ and

$$\gamma_1(s) := \exp(\sigma r)L\left(s\exp(\sigma r)\right)+\varphi\left(s\exp(\sigma r)\right) \tag{5.9a}$$

$$b_1(s) := \frac{1}{\varphi\left(s\exp(\sigma r)\right)} \tag{5.9b}$$

Step $i \geq 2$: Based on the knowledge of the functions $\mu_j \in C^\infty(\Re^j;(0,+\infty))$, $\gamma_j \in C^\infty(\Re^+;(0,+\infty))$, $b_j \in C^\infty(\Re^+;(0,+\infty))$ ($j = 1,...,i-1$) we will define the functions $\mu_i \in C^\infty(\Re^i;(0,+\infty))$, $\gamma_i \in C^\infty(\Re^+;(0,+\infty))$, $b_i \in C^\infty(\Re^+;(0,+\infty))$ such that the following inequality holds for all $(\xi_1,...,\xi_i)' \in \Re^i$:

$$\begin{aligned}
&-\xi_1^2 b_1(s)\mu_1(\xi_1) + s|\xi_1|\gamma_1(s) - \sum_{j=2}^{i}(\xi_j - k_{j-1}(\xi_1,...,\xi_{j-1}))^2 b_j(s)\mu_j(\xi_1,...,\xi_j) \\
&+ s\sum_{j=2}^{i}\left|\xi_j - k_{j-1}(\xi_1,...,\xi_{j-1})\right|\gamma_j(s) + s\sum_{j=2}^{i}\left|\xi_j - k_{j-1}(\xi_1,...,\xi_{j-1})\right|\left(\sum_{k=1}^{j-1}\gamma_k(s)\right)\delta_{j-1}(\xi_1,...,\xi_{j-1}) \\
&\leq -(n+1-i)\sigma\left(\xi_1^2 + \sum_{j=2}^{i}\left|\xi_j - k_{j-1}(\xi_1,...,\xi_{j-1})\right|^2\right)
\end{aligned} \tag{5.10}$$

where

$$k_1(\xi_1) := -\mu_1(\xi_1)\,\xi_1 \tag{5.11a}$$

$$k_j(\xi_1,...,\xi_j) := -\mu_j(\xi_1,...,\xi_j)\left(\xi_j - k_{j-1}(\xi_1,...,\xi_{j-1})\right),\ j = 2,...,i-1 \tag{5.11b}$$

$$s = |\xi_1| + \sum_{j=2}^{i}\left|\xi_j - k_{j-1}(\xi_1,...,\xi_{j-1})\right| \tag{5.12}$$

$$\delta_j(\xi_1,...,\xi_j) := \left|\nabla k_j(\xi_1,...,\xi_j)\right|\left(1 + \mu_1(\xi_1) + ... + \mu_j(\xi_1,...,\xi_i)\right),\ j = 1,...,i-1 \tag{5.13}$$

and

$$\gamma_i(s) := i\exp(\sigma r)L\left(is\exp(\sigma r)B_i\left(is\exp(\sigma r)\right)\right)B_i\left(is\exp(\sigma r)\right) + \varphi\left(is\exp(\sigma r)B_i\left(is\exp(\sigma r)\right)\right) \tag{5.14a}$$

$$b_i(s) := \frac{1}{\varphi\left(is\exp(\sigma r)B_i\left(is\exp(\sigma r)\right)\right)} \tag{5.14b}$$

where $B_i \in C^\infty(\Re^+;(0,+\infty))$ is a non-decreasing function that satisfies:

$$B_i(s) \geq \max\left\{1 + \sum_{j=1}^{i-1}\mu_j(\xi_1,...,\xi_j) : |\xi_1| + \sum_{j=2}^{i}\left|\xi_j - k_{j-1}(\xi_1,...,\xi_{j-1})\right| \leq s\right\}\ \text{, for all } s \geq 0 \tag{5.14c}$$

Indeed, by the previous step we get the existence of functions $\mu_j \in C^\infty(\Re^j;(0,+\infty))$, $\gamma_j \in C^\infty(\Re^+;(0,+\infty))$, $b_j \in C^\infty(\Re^+;(0,+\infty))$ ($j = 1,...,i-1$) such that the following inequalities hold for all $(\xi_1,...,\xi_{i-1})' \in \Re^{i-1}$

$$-\xi_1^2 b_1(s')\mu_1(\xi_1) + s'|\xi_1|\gamma_1(s') \leq -n\sigma\xi_1^2\ \text{, for the case } i = 2 \tag{5.15a}$$

and for the case $i > 2$:

$$-\xi_1^2 b_1(s')\mu_1(\xi_1)+s'|\xi_1|\gamma_1(s')-\sum_{j=2}^{i-1}(\xi_j-k_{j-1}(\xi_1,...,\xi_{j-1}))^2 b_j(s')\mu_j(\xi_1,...,\xi_j)$$
$$+s'\sum_{j=2}^{i-1}\left|\xi_j-k_{j-1}(\xi_1,...,\xi_{j-1})\right|\gamma_j(s')+s'\sum_{j=2}^{i-1}\left|\xi_j-k_{j-1}(\xi_1,...,\xi_{j-1})\right|\left(\sum_{k=1}^{j-1}\gamma_k(s')\right)\delta_{j-1}(\xi_1,...,\xi_{j-1}) \tag{5.15b}$$
$$\leq -(n+2-i)\sigma\left(\xi_1^2+\sum_{j=2}^{i-1}\left|\xi_j-k_{j-1}(\xi_1,...,\xi_{j-1})\right|^2\right)$$

where

$$s'=|\xi_1|\text{, for the case } i=2 \tag{5.16a}$$

$$s'=|\xi_1|+\sum_{j=2}^{i-1}\left|\xi_j-k_{j-1}(\xi_1,...,\xi_{j-1})\right|\text{, for the case } i>2 \tag{5.16b}$$

Let functions $\rho_j\in C^\infty(\Re^+;(0,+\infty))$ ( $j=1,...,i$ ) such that the following inequalities hold for all $s\geq s'\geq 0$ :

$$b_j(s)-b_j(s')+s\gamma_j(s)-s'\gamma_j(s')\leq (s-s')\rho_j(s) \tag{5.17}$$

By virtue of inequalities (5.15), (5.17) and definitions (5.12), (5.16), we obtain for all $(\xi_1,...,\xi_i)'\in\Re^i$ :

$$-\xi_1^2 b_1(s)\mu_1(\xi_1)+s|\xi_1|\gamma_1(s)-\sum_{j=2}^{i}(\xi_j-k_{j-1}(\xi_1,...,\xi_{j-1}))^2 b_j(s)\mu_j(\xi_1,...,\xi_j)$$
$$+s\sum_{j=2}^{i}\left|\xi_j-k_{j-1}(\xi_1,...,\xi_{j-1})\right|\gamma_j(s)+s\sum_{j=2}^{i}\left|\xi_j-k_{j-1}(\xi_1,...,\xi_{j-1})\right|\left(\sum_{k=1}^{j-1}\gamma_k(s)\right)\delta_{j-1}(\xi_1,...,\xi_{j-1}) \tag{5.18a}$$
$$\leq -n\sigma\xi_1^2-(\xi_2-k_1(\xi_1))^2 b_2(s)\mu_2(\xi_1,\xi_2)+s|\xi_2-k_1(\xi_1)|\gamma_2(s)$$
$$+s|\xi_2-k_1(\xi_1)|\gamma_1(s)\delta_1(\xi_1)+|\xi_1|\rho_1(s)|\xi_2-k_1(\xi_1)|+\rho_1(s)|\xi_2-k_1(\xi_1)|\mu_1(\xi_1)\xi_1^2$$

for the case $i=2$ and

$$-\xi_1^2 b_1(s)\mu_1(\xi_1)+s|\xi_1|\gamma_1(s)-\sum_{j=2}^{i}(\xi_j-k_{j-1}(\xi_1,...,\xi_{j-1}))^2 b_j(s)\mu_j(\xi_1,...,\xi_j)$$
$$+s\sum_{j=2}^{i}\left|\xi_j-k_{j-1}(\xi_1,...,\xi_{j-1})\right|\gamma_j(s)+s\sum_{j=2}^{i}\left|\xi_j-k_{j-1}(\xi_1,...,\xi_{j-1})\right|\left(\sum_{k=1}^{j-1}\gamma_k(s)\right)\delta_{j-1}(\xi_1,...,\xi_{j-1})$$
$$\leq -(n+2-i)\sigma\left(\xi_1^2+\sum_{j=2}^{i-1}\left|\xi_j-k_{j-1}(\xi_1,...,\xi_{j-1})\right|^2\right)-(\xi_i-k_{i-1}(\xi_1,...,\xi_{i-1}))^2 b_i(s)\mu_i(\xi_1,...,\xi_i)$$
$$+|\xi_1|\rho_1(s)|\xi_i-k_{i-1}(\xi_1,...,\xi_{i-1})|+s|\xi_i-k_{i-1}(\xi_1,...,\xi_{i-1})|\left(\sum_{k=1}^{i-1}\gamma_k(s)\right)\delta_{i-1}(\xi_1,...,\xi_{j-1})$$
$$+|\xi_i-k_{i-1}(\xi_1,...,\xi_{i-1})|\sum_{j=2}^{i-1}\left|\xi_j-k_{j-1}(\xi_1,...,\xi_{j-1})\right|\left(\sum_{k=1}^{j-1}\rho_k(s)\right)\delta_{j-1}(\xi_1,...,\xi_{j-1})$$
$$+|\xi_i-k_{i-1}(\xi_1,...,\xi_{i-1})|\sum_{j=2}^{i-1}\left|\xi_j-k_{j-1}(\xi_1,...,\xi_{j-1})\right|\rho_j(s)+s|\xi_i-k_{i-1}(\xi_1,...,\xi_{i-1})|\gamma_i(s)$$
$$+\xi_1^2\rho_1(s)\mu_1(\xi_1)|\xi_i-k_{i-1}(\xi_1,...,\xi_{i-1})|+|\xi_i-k_{i-1}(\xi_1,...,\xi_{i-1})|\sum_{j=2}^{i-1}(\xi_j-k_{j-1}(\xi_1,...,\xi_{j-1}))^2\rho_j(s)\mu_j(\xi_1,...,\xi_j)$$
$$\tag{5.18b}$$

for the case $i>2$. Completing the squares in the right-hand sides of inequalities (5.18a,b) and using definition (5.12), we get:

$$
\begin{aligned}
&s|\xi_2-k_1(\xi_1)|\gamma_2(s)+s|\xi_2-k_1(\xi_1)|\gamma_1(s)\delta_1(\xi_1)+|\xi_1|\rho_1(s)|\xi_2-k_1(\xi_1)|+\rho_1(s)|\xi_2-k_1(\xi_1)|\mu_1(\xi_1)\xi_1^2\\
&\leq|\xi_2-k_1(\xi_1)|^2\gamma_2(s)+|\xi_2-k_1(\xi_1)|^2\gamma_1(s)\delta_1(\xi_1)+\sigma\xi_1^2+\frac{1}{\sigma}|\xi_2-k_1(\xi_1)|^2\gamma_2^2(s)\\
&+\frac{1}{\sigma}|\xi_2-k_1(\xi_1)|^2\gamma_1^2(s)\delta_1^2(\xi_1)+\frac{1}{\sigma}\rho_1^2(s)|\xi_2-k_1(\xi_1)|^2+\frac{1}{\sigma}\rho_1^2(s)|\xi_2-k_1(\xi_1)|^2\mu_1^2(\xi_1)\xi_1^2
\end{aligned}
\tag{5.19a}
$$

for the case $i=2$ and

$$
\begin{aligned}
&|\xi_1|\rho_1(s)|\xi_i-k_{i-1}(\xi_1,...,\xi_{i-1})|+s|\xi_i-k_{i-1}(\xi_1,...,\xi_{i-1})|\left(\sum_{k=1}^{i-1}\gamma_k(s)\right)\delta_{i-1}(\xi_1,...,\xi_{j-1})\\
&+|\xi_i-k_{i-1}(\xi_1,...,\xi_{i-1})|\sum_{j=2}^{i-1}|\xi_j-k_{j-1}(\xi_1,...,\xi_{j-1})|\left(\sum_{k=1}^{j-1}\rho_k(s)\right)\delta_{j-1}(\xi_1,...,\xi_{j-1})\\
&+|\xi_i-k_{i-1}(\xi_1,...,\xi_{i-1})|\sum_{j=2}^{i-1}|\xi_j-k_{j-1}(\xi_1,...,\xi_{j-1})|\rho_j(s)+s|\xi_i-k_{i-1}(\xi_1,...,\xi_{i-1})|\gamma_i(s)\\
&+\xi_1^2\rho_1(s)\mu_1(\xi_1)|\xi_i-k_{i-1}(\xi_1,...,\xi_{i-1})|+|\xi_i-k_{i-1}(\xi_1,...,\xi_{i-1})|\sum_{j=2}^{i-1}(\xi_j-k_{j-1}(\xi_1,...,\xi_{j-1}))^2\rho_j(s)\mu_j(\xi_1,...,\xi_j)\\
&\leq|\xi_i-k_{i-1}(\xi_1,...,\xi_{i-1})|^2\left(\sum_{k=1}^{i-1}\gamma_k(s)\right)\delta_{i-1}(\xi_1,...,\xi_{j-1})\\
&+\frac{5(i-1)}{4\sigma}|\xi_i-k_{i-1}(\xi_1,...,\xi_{i-1})|^2\left(\sum_{k=1}^{i-1}\gamma_k(s)\right)^2\delta_{i-1}^2(\xi_1,...,\xi_{j-1})+\\
&+\frac{5}{4\sigma}|\xi_i-k_{i-1}(\xi_1,...,\xi_{i-1})|^2\sum_{j=2}^{i-1}\left(\sum_{k=1}^{j-1}\rho_k(s)\right)^2\delta_{j-1}^2(\xi_1,...,\xi_{j-1})\\
&+\frac{5}{4\sigma}\left(\sum_{j=1}^{i-1}\rho_j^2(s)\right)|\xi_i-k_{i-1}(\xi_1,...,\xi_{i-1})|^2+\sigma\left(|\xi_1|^2+\sum_{j=1}^{i-1}|\xi_j-k_{j-1}(\xi_1,...,\xi_{j-1})|^2\right)\\
&+|\xi_i-k_{i-1}(\xi_1,...,\xi_{i-1})|^2\gamma_i(s)+\frac{5(i-1)}{4\sigma}|\xi_i-k_{i-1}(\xi_1,...,\xi_{i-1})|^2\gamma_i^2(s)\\
&+\frac{5}{4\sigma}\xi_1^2\rho_1^2(s)\mu_1^2(\xi_1)|\xi_i-k_{i-1}(\xi_1,...,\xi_{i-1})|^2\\
&+\frac{5}{4\sigma}|\xi_i-k_{i-1}(\xi_1,...,\xi_{i-1})|^2\sum_{j=2}^{i-1}(\xi_j-k_{j-1}(\xi_1,...,\xi_{j-1}))^2\rho_j^2(s)\mu_j^2(\xi_1,...,\xi_j)
\end{aligned}
\tag{5.19b}
$$

for the case $i>2$. It follows from (5.18), (5.19) that the selection:

$$
\mu_2(\xi_1,\xi_2)=b_2^{-1}(p)\left[(n-1)\sigma+\gamma_2(p)+\gamma_1(p)\delta_1(\xi_1)+\frac{1}{\sigma}\gamma_2^2(p)+\frac{1}{\sigma}\gamma_1^2(p)\delta_1^2(\xi_1)+\frac{1}{\sigma}\left(1+\mu_1^2(\xi_1)\xi_1^2\right)\rho_1^2(p)\right]
\tag{5.20a}
$$

for the case $i=2$, where

$$
p:=1+\xi_1^2+|\xi_2-k_1(\xi_1)|^2
\tag{5.20b}
$$

and

$$\mu_i(\xi_1,...,\xi_i):=$$

$$b_i^{-1}(p)\left[(n+1-i)\sigma+\frac{5(i-1)}{4\sigma}\gamma_i^2(p)+\gamma_i(p)+\left(\sum_{k=1}^{i-1}\gamma_k(p)\right)\delta_{i-1}(\xi_1,...,\xi_{j-1})+\frac{5}{4\sigma}\sum_{j=2}^{i-1}(\xi_j-k_{j-1}(\xi_1,...,\xi_{j-1}))^2\rho_j^2(p)\mu_j^2(\xi_1,...,\xi_j)\right]$$

$$+\frac{5}{4\sigma}b_i^{-1}(p)\left[\frac{i-1}{\sigma}\left(\sum_{k=1}^{i-1}\gamma_k(p)\right)^2\delta_{i-1}^2(\xi_1,...,\xi_{j-1})+\left(\sum_{j=1}^{i-1}\rho_j^2(p)\right)+\sum_{j=2}^{i-1}\left(\sum_{k=1}^{j-1}\rho_k(p)\right)^2\delta_{j-1}^2(\xi_1,...,\xi_{j-1})+\xi_1^2\rho_1^2(p)\mu_1^2(\xi_1)\right] \tag{5.21a}$$

for the case $i>2$, where

$$p:=\frac{i}{2}+\left(\xi_1^2+\sum_{j=2}^{i}\left|\xi_j-k_{j-1}(\xi_1,...,\xi_{j-1})\right|^2\right) \tag{5.21b}$$

guarantees inequality (5.10).

Having performed $n$ steps of the above algorithm, we have defined functions $\mu_i\in C^\infty(\Re^i;(0,+\infty))$, $\gamma_i\in C^\infty(\Re^+;(0,+\infty))$, $b_i\in C^\infty(\Re^+;(0,+\infty))$, $k_i\in C^\infty(\Re^i;\Re)$ ($i=1,...,n$) with

$$k_1(\xi_1):=-\mu_1(\xi_1)\,\xi_1 \tag{5.22a}$$

$$k_j(\xi_1,...,\xi_j):=-\mu_j(\xi_1,...,\xi_j)\left(\xi_j-k_{j-1}(\xi_1,...,\xi_{j-1})\right),\ j=2,...,n \tag{5.22b}$$

$$\gamma_1(s):=\exp(\sigma\,r)L\left(s\exp(\sigma\,r)\right)+\varphi\left(s\exp(\sigma\,r)\right) \tag{5.23a}$$

$$\gamma_i(s):=i\exp(\sigma\,r)L\left(is\exp(\sigma\,r)B_i\left(is\exp(\sigma\,r)\right)\right)B_i\left(is\exp(\sigma\,r)\right)+\varphi\left(is\exp(\sigma\,r)B_i\left(is\exp(\sigma\,r)\right)\right),\ i=2,...,n \tag{5.23b}$$

$$b_i(s):=\frac{1}{\varphi\left(is\exp(\sigma\,r)B_i\left(is\exp(\sigma\,r)\right)\right)} \tag{5.23c}$$

where $B_i\in C^\infty(\Re^+;(0,+\infty))$ ($i=2,...,n$) are non-decreasing functions that satisfy

$$B_i(s)\geq\max\left\{1+\sum_{j=1}^{i-1}\mu_j(\xi_1,...,\xi_j):\left|\xi_1\right|+\sum_{j=2}^{i}\left|\xi_j-k_{j-1}(\xi_1,...,\xi_{j-1})\right|\leq s\right\}\ \text{, for all }s\geq 0 \tag{5.23d}$$

and such that

$$\begin{aligned}&-\xi_1^2b_1(s)\mu_1(\xi_1)+s\left|\xi_1\right|\gamma_1(s)-\sum_{j=2}^{n}(\xi_j-k_{j-1}(\xi_1,...,\xi_{j-1}))^2b_j(s)\mu_j(\xi_1,...,\xi_j)\\&+s\sum_{j=2}^{n}\left|\xi_j-k_{j-1}(\xi_1,...,\xi_{j-1})\right|\gamma_j(s)+s\sum_{j=2}^{n}\left|\xi_j-k_{j-1}(\xi_1,...,\xi_{j-1})\right|\left(\sum_{k=1}^{j-1}\gamma_k(s)\right)\delta_{j-1}(\xi_1,...,\xi_{j-1})\\&\leq-\sigma\left(\xi_1^2+\sum_{j=2}^{n}\left|\xi_j-k_{j-1}(\xi_1,...,\xi_{j-1})\right|^2\right)\end{aligned} \tag{5.24}$$

where

$$s=\left|\xi_1\right|+\sum_{j=2}^{n}\left|\xi_j-k_{j-1}(\xi_1,...,\xi_{j-1})\right| \tag{5.25a}$$

$$\delta_i(\xi_1,...,\xi_i) := \left|\nabla k_i(\xi_1,...,\xi_i)\right|\left(1+\mu_1(\xi_1)+...+\mu_i(\xi_1,...,\xi_i)\right),\ i=1,...,n-1 \quad (5.25b)$$

By virtue of Lemma 5.3 it follows that the functional $V$ defined by (5.4) satisfies:

$$V^0(x;v) \le -2\sigma V(x)$$

$$\text{for all } (t,x,u,d) \in \Re^+ \times C^0([-r,0];\Re^n) \times \Re \times D \text{ with } V(x) > x_1^2(0) + \sum_{j=2}^{n} \left|x_j(0) - k_{j-1}(x_1(0),...,x_{j-1}(0))\right|^2$$

$$\text{and } v = (f_1(t,d,x_1) + g_1(t,d,x_1)x_2(0),...,f_n(t,d,x) + g_n(t,d,x)u)' \in \Re^n \quad (5.26)$$

$$V^0(x;v) \le \max\left\{-2\sigma V(x), 2A(t,d,x,u)\right\}$$

$$\text{for all } (t,x,u,d) \in \Re^+ \times C^0([-r,0];\Re^n) \times \Re \times D \text{ with } V(x) = x_1^2(0) + \sum_{j=2}^{n} \left|x_j(0) - k_{j-1}(x_1(0),...,x_{j-1}(0))\right|^2$$

$$\text{and } v = (f_1(t,d,x_1) + g_1(t,d,x_1)x_2(0),...,f_n(t,d,x) + g_n(t,d,x)u)' \in \Re^n \quad (5.27)$$

where

$$\begin{aligned}
&A(t,d,x,u) := x_1(0)\left(f_1(t,d,x_1) + g_1(t,d,x_1)x_2(0)\right) \\
&+\sum_{j=2}^{n-1}\left(x_j(0) - k_{j-1}(x_1(0),...,x_{j-1}(0))\right)\left(f_i(t,d,x_1,...,x_i) + g_i(t,d,x_1,...,x_i)x_{i+1}(0)\right) \\
&-\sum_{j=2}^{n-1}\left(x_j(0) - k_{j-1}(x_1(0),...,x_{j-1}(0))\right)\left(\sum_{l=1}^{j-1}\frac{\partial k_{j-1}}{\partial \xi_l}(x_1(0),...,x_{j-1}(0))\left(f_l(t,d,x_1,...,x_l) + g_l(t,d,x_1,...,x_l)x_{l+1}(0)\right)\right) \quad (5.28) \\
&+\left(x_n(0) - k_{n-1}(x_1(0),...,x_{n-1}(0))\right)\left(f_n(t,d,x) + g_n(t,d,x)u\right) \\
&-\left(x_n(0) - k_{n-1}(x_1(0),...,x_{n-1}(0))\right)\sum_{l=1}^{n-1}\frac{\partial k_{n-1}}{\partial \xi_l}(x_1(0),...,x_{n-1}(0))\left(f_l(t,d,x_1,...,x_l) + g_l(t,d,x_1,...,x_l)x_{l+1}(0)\right)
\end{aligned}$$

We next show that $A(t,d,x,u) \le -\sigma\left(x_1^2(0) + \sum_{j=2}^{n}\left|x_j(0) - k_{j-1}(x_1(0),...,x_{j-1}(0))\right|^2\right)$ for all $(t,x,d) \in \Re^+ \times C^0([-r,0];\Re^n) \times D$ with $u = k_n(x_1(0),...,x_n(0))$ and $V(x) = x_1^2(0) + \sum_{j=2}^{n}\left|x_j(0) - k_{j-1}(x_1(0),...,x_{j-1}(0))\right|^2$. Indeed, notice that by virtue of definitions (5.11), (5.12) for every $i \ge 2$ and $(\xi_1,...,\xi_i)' \in \Re^i$ it holds that

$$\sum_{j=1}^{i}\left|\xi_j\right| \le \sum_{j=1}^{i}\left|\xi_j - k_{j-1}(\xi_1,...,\xi_{j-1})\right| + \sum_{j=2}^{i}\mu_{j-1}(\xi_1,...,\xi_{j-1})\left|\xi_{j-1} - k_{j-2}(\xi_1,...,\xi_{j-2})\right| \quad (5.29)$$

Notice that here we have used the convention $k_0 \equiv 0$. Inequality (5.29) implies:

$$\sum_{j=1}^{i}\left|\xi_j\right| \le \left|\xi_i - k_{i-1}(\xi_1,...,\xi_{i-1})\right| + \sum_{j=1}^{i-1}\left(1+\mu_j(\xi_1,...,\xi_{j-1})\right)\left|\xi_j - k_{j-1}(\xi_1,...,\xi_{j-1})\right| \quad (5.30)$$

Consequently, inequalities (5.1), (5.7) in conjunction with (5.30) and (5.23d) imply for all $i \ge 2$, $(t,(x_1,...,x_i),d) \in \Re^+ \times C^0([-r,0];\Re^i) \times D$ :

$$\left|f_i(t,d,x_1,...,x_i)\right| \le \max_{\theta\in[-r,0]} L\left(B_i\left(P_i x(\theta)\right)P_i x(\theta)\right)B_i\left(P_i x(\theta)\right)P_i x(\theta) \quad (5.31a)$$

$$\frac{1}{\varphi\big(B_i\left(P_i x(\theta)\right)P_i x(\theta)\big)} \le g_i(t,d,x_1,...,x_i) \le \varphi\big(B_i\left(P_i x(\theta)\right)P_i x(\theta)\big) \tag{5.31b}$$

where

$$P_i x(\theta) := \left|x_1(\theta)\right| + \sum_{j=2}^{i} \left|x_j(\theta) - k_{j-1}(x_1(\theta),...,x_{j-1}(\theta))\right| \tag{5.32}$$

The above inequalities in conjunction with the fact $\max_{\theta\in[-r,0]} \exp(2\sigma\theta)\left( x_1^2(\theta) + \sum_{j=2}^{n}\left|x_j(\theta) - k_{j-1}(x_1(\theta),...,x_{j-1}(\theta))\right|^2 \right) = x_1^2(0) + \sum_{j=2}^{n}\left|x_j(0) - k_{j-1}(x_1(0),...,x_{j-1}(0))\right|^2$ (which implies $P_i x(\theta) \le i\exp(\sigma r)\left( \left|x_1(0)\right| + \sum_{j=2}^{n}\left|x_j(0) - k_{j-1}(x_1(0),...,x_{j-1}(0))\right| \right)$ for all $\theta\in[-r,0]$) give for all $i \ge 2$, $(t,(x_1,...,x_i),d) \in \Re^+ \times C^0([-r,0];\Re^i)\times D$ and $V(x) = x_1^2(0) + \sum_{j=2}^{n}\left|x_j(0) - k_{j-1}(x_1(0),...,x_{j-1}(0))\right|^2$ :

$$\left|f_i(t,d,x_1,...,x_i)\right| \le L\big(is\exp(\sigma r)B_i\left(is\exp(\sigma r)\right)\big)B_i\left(is\exp(\sigma r)\right)is\exp(\sigma r) \tag{5.31a}$$

$$\frac{1}{\varphi\big(is\exp(\sigma r)B_i\left(is\exp(\sigma r)\right)\big)} \le g_i(t,d,x_1,...,x_i) \le \varphi\big(is\exp(\sigma r)B_i\left(is\exp(\sigma r)\right)\big) \tag{5.31b}$$

where

$$s := \left|x_1(0)\right| + \sum_{j=2}^{n}\left|x_j(0) - k_{j-1}(x_1(0),...,x_{j-1}(0))\right| \tag{5.32}$$

Using (5.31a,b), (5.32) and performing a similar analysis for the case $i=1$, we obtain by virtue of definitions (5.23a,b,c):

$$\left|f_i(t,d,x_1,...,x_i)\right| + \left|g_i(t,d,x_1,...,x_i)\right|\left|x_{i+1}(0) - k_i(x_1(0),...,x_i(0))\right| \le s\gamma_i(s)\,,\ i=1,...,n-1 \tag{5.33a}$$

$$\left|f_n(t,d,x_1,...,x_n)\right| \le s\gamma_n(s) \tag{5.33b}$$

$$b_i(s) \le g_i(t,d,x_1,...,x_i) \le \gamma_i(s)\,,\ i=1,...,n \tag{5.33c}$$

Definition (5.28) in conjunction with definitions (5.22) gives for $u = k_n(x_1(0),...,x_n(0))$ :

$$
\begin{aligned}
A(t,d,x,u) \le & -g_1(t,d,x_1)\mu_1(x_1(0))x_1^2(0) + |x_1(0)|\left| f_1(t,d,x_1) + g_1(t,d,x_1)\left(x_2(0) - k_1(x_1(0))\right)\right| \\
& - \sum_{j=2}^{n} g_i(t,d,x_1,...,x_i)\mu_j(x_1(0),...,x_j(0))\left(x_j(0) - k_{j-1}(x_1(0),...,x_{j-1}(0))\right)^2 \\
& + \sum_{j=2}^{n-1} \left|x_j(0) - k_{j-1}(x_1(0),...,x_{j-1}(0))\right| \left| f_i(t,d,x_1,...,x_i) + g_i(t,d,x_1,...,x_i)(x_{i+1}(0) - k_j(x_1(0),...,x_j(0)))\right| \\
& + \sum_{j=2}^{n-1} \left|\nabla k_{j-1}(x_1(0),...,x_{j-1}(0))\right| \left|x_j(0) - k_{j-1}(x_1(0),...,x_{j-1}(0))\right| \left| \sum_{l=1}^{j-1} f_l(t,d,x_1,...,x_l) + g_l(t,d,x_1,...,x_l)(x_{l+1}(0) - k_l(x_1(0),...,x_l(0)))\right| \\
& + \sum_{j=2}^{n-1} \left|\nabla k_{j-1}(x_1(0),...,x_{j-1}(0))\right| \left|x_j(0) - k_{j-1}(x_1(0),...,x_{j-1}(0))\right| \left| \sum_{l=1}^{j-1} \left|x_l(0) - k_{l-1}(x_1(0),...,x_{l-1}(0))\right| \mu_l(x_1(0),...,x_l(0)) \right| g_l(t,d,x_1,...,x_l)\Big| \\
& + \left|x_n(0) - k_{n-1}(x_1(0),...,x_{n-1}(0))\right| \left| f_n(t,d,x)\right| \\
& + \left|\nabla k_{n-1}(x_1(0),...,x_{n-1}(0))\right| \left|x_n(0) - k_{n-1}(x_1(0),...,x_{n-1}(0))\right| \left| \sum_{l=1}^{n-1} f_l(t,d,x_1,...,x_l) + g_l(t,d,x_1,...,x_l)(x_{l+1}(0) - k_l(x_1(0),...,x_l(0)))\right| \\
& + \left|\nabla k_{n-1}(x_1(0),...,x_{n-1}(0))\right| \left|x_n(0) - k_{n-1}(x_1(0),...,x_{n-1}(0))\right| \left| \sum_{l=1}^{n-1} \left|x_l(0) - k_{l-1}(x_1(0),...,x_{l-1}(0))\right| \mu_l(x_1(0),...,x_l(0)) \right| g_l(t,d,x_1,...,x_l)\Big|
\end{aligned}
\tag{5.34}
$$

Combining (5.34) with (5.33a,b,c) we obtain for $u = k_n(x_1(0),..., x_n(0))$ :

$$
\begin{aligned}
A(t,d,x,u) \le & -b_1(s)\mu_1(x_1(0))x_1^2(0) + |x_1(0)| s\gamma_1(s) \\
& - \sum_{j=2}^{n} b_j(s)\mu_j(x_1(0),...,x_j(0))\left(x_j(0) - k_{j-1}(x_1(0),...,x_{j-1}(0))\right)^2 \\
& + s\sum_{j=2}^{n} \left|x_j(0) - k_{j-1}(x_1(0),...,x_{j-1}(0))\right| \gamma_j(s) \\
& + s\sum_{j=2}^{n} \left|\nabla k_{j-1}(x_1(0),...,x_{j-1}(0))\right| \left|x_j(0) - k_{j-1}(x_1(0),...,x_{j-1}(0))\right| \sum_{l=1}^{j-1} \left(1 + \mu_l(x_1(0),...,x_l(0))\right)\gamma_l(s)
\end{aligned}
\tag{5.35}
$$

The above inequality in conjunction with definition (5.25b) implies for $u = k_n(x_1(0),..., x_n(0))$ :

$$
\begin{aligned}
A(t,d,x,u) \le & -b_1(s)\mu_1(x_1(0))x_1^2(0) + |x_1(0)| s\gamma_1(s) \\
& - \sum_{j=2}^{n} b_j(s)\mu_j(x_1(0),...,x_j(0))\left(x_j(0) - k_{j-1}(x_1(0),...,x_{j-1}(0))\right)^2 \\
& + s\sum_{j=2}^{n} \left|x_j(0) - k_{j-1}(x_1(0),...,x_{j-1}(0))\right| \gamma_j(s) \\
& + s\sum_{j=2}^{n} \delta_{j-1}(x_1(0),...,x_{j-1}(0)) \left|x_j(0) - k_{j-1}(x_1(0),...,x_{j-1}(0))\right| \sum_{l=1}^{j-1} \gamma_l(s)
\end{aligned}
\tag{5.36}
$$

Clearly, inequality (5.36) in conjunction with (5.24) show that $A(t,d,x,u) \le -\sigma\left( x_1^2(0) + \sum_{j=2}^{n} \left|x_j(0) - k_{j-1}(x_1(0),...,x_{j-1}(0))\right|^2 \right)$ for all $(t,x,d) \in \Re^+ \times C^0([-r,0];\Re^n) \times D$ with $u = k_n(x_1(0),..., x_n(0))$ and $V(x) = x_1^2(0) + \sum_{j=2}^{n} \left|x_j(0) - k_{j-1}(x_1(0),..., x_{j-1}(0))\right|^2$. By virtue of (5.26), (5.27) we obtain:

$$V^0(x;v) \le -2\sigma V(x)$$
for all $(t,x,d) \in \Re^+ \times C^0([-r,0];\Re^n) \times D$ with
$$v = (f_1(t,d,x_1) + g_1(t,d,x_1)x_2(0),\ldots, f_n(t,d,x) + g_n(t,d,x)k_n(x(0)))' \in \Re^n \tag{5.37}$$

The reader should notice that there exist functions $a_1, a_2 \in K_\infty$ such that

$$a_1(|\xi|) \le \xi_1^2 + \sum_{j=2}^{n} \left|\xi_j - k_{j-1}(\xi_1,\ldots,\xi_{j-1})\right|^2 \le a_2(|\xi|),\ \forall \xi = (\xi_1,\ldots,\xi_n)' \in \Re^n \tag{5.38}$$

Consequently, definition (5.4) implies

$$\exp(-2\sigma r)\, a_1(\|x\|_r) \le V(x) \le a_2(\|x\|_r),\ \forall x \in C^0([-r,0];\Re^n) \tag{5.39}$$

It follows from inequalities (5.37), (5.39) and Proposition 4.6 that the closed-loop system (1.4) with $u(t) = k_n(x(t))$ is URGAS.

Finally, we show that $V$ as defined by (5.4) is a SRCLF, which satisfies the "small-control" property. Clearly, definition (5.28) in conjunction with the fact that $A(t,d,x,u) \le -\sigma\left( x_1^2(0) + \sum_{j=2}^{n} \left|x_j(0) - k_{j-1}(x_1(0),\ldots,x_{j-1}(0))\right|^2 \right)$ for all $(t,x,d) \in \Re^+ \times C^0([-r,0];\Re^n) \times D$ with $u = k_n(x_1(0),\ldots,x_n(0))$ and $V(x) = x_1^2(0) + \sum_{j=2}^{n} \left|x_j(0) - k_{j-1}(x_1(0),\ldots,x_{j-1}(0))\right|^2$, implies for all $(t,x,u,d) \in \Re^+ \times C^0([-r,0];\Re^n) \times \Re \times D$ with $V(x) = x_1^2(0) + \sum_{j=2}^{n} \left|x_j(0) - k_{j-1}(x_1(0),\ldots,x_{j-1}(0))\right|^2$:

$$\begin{aligned} A(t,d,x,u) \le & -\sigma\left( x_1^2(0) + \sum_{j=2}^{n} \left|x_j(0) - k_{j-1}(x_1(0),\ldots,x_{j-1}(0))\right|^2 \right) \\ & + \left|x_n(0) - k_{n-1}(x_1(0),\ldots,x_{n-1}(0))\right| \left|g_n(t,d,x)\right| \left|u - k_n(x_1(0),\ldots,x_n(0))\right| \end{aligned} \tag{5.40}$$

By virtue of (5.33c), (5.26), (5.27) and (5.40) we obtain:

$$V^0(x;v) \le -2\sigma V(x) + 2\left|x_n(0) - k_{n-1}(x_1(0),\ldots,x_{n-1}(0))\right| \gamma_n\left( |x_1(0)| + \sum_{j=2}^{n} \left|x_j(0) - k_{j-1}(x_1(0),\ldots,x_{j-1}(0))\right| \right) \left|u - k_n(x_1(0),\ldots,x_n(0))\right|$$
for all $(t,x,u,d) \in \Re^+ \times C^0([-r,0];\Re^n) \times \Re \times D$ with
$$v = (f_1(t,d,x_1) + g_1(t,d,x_1)x_2(0),\ldots, f_n(t,d,x) + g_n(t,d,x)u)' \in \Re^n \tag{5.41}$$

Define:

$$\Phi(t,x) := x(0),\ \rho(w) := 2\sigma w \text{ and}$$
$$\Psi(t,z,u) := 2\left|z_n - k_{n-1}(z_1,\ldots,z_{n-1})\right| \gamma_n\left( |z_1| + \sum_{j=2}^{n} \left|z_j - k_{j-1}(z_1,\ldots,z_{j-1})\right| \right) \left|u - k_n(z_1,\ldots,z_n)\right|$$

The above definitions in conjunction with inequalities (5.39) and (5.41) guarantee that inequalities (4.7), (4.10) and (4.11) hold, for $V$ as defined by (5.4). Consequently, $V$ as defined by (5.4) is a SRCLF, which satisfies the "small-control" property.

The proof is complete. ⊲

**Example 5.4:** Consider the control system:

$$\begin{aligned}
&\dot{x}_1(t) = d_1(t)\int_{t-r}^{t} x_1^2(\theta)d\theta + x_2(t) \\
&\dot{x}_2(t) = d_2(t)\|T_r(t)x_2\|_r + u(t) \\
&(x_1(t), x_2(t)) \in \Re^2 \ , (d_1(t), d_2(t)) \in [-1,1]^2 \ , u(t) \in \Re
\end{aligned} \tag{5.42}$$

Clearly, system (5.42) is a control system described by RFDEs, which satisfies the hypotheses of Theorem 5.1. More specifically, inequality (5.1) holds with $b_1 = b_2 = \varphi \equiv 1$. In order to design a delay free stabilizing feedback for (5.42) we follow the algorithm in the proof of Theorem 5.1. Notice that inequality (5.7) holds with $L(w) = 1 + r\,w$. Let $\sigma > 0$ be given.

Step $i = 1$: We define:

$$\mu_1(\xi_1) := \exp(\sigma\, r)\left(1 + r\,(1+\xi_1^2)\exp(\sigma\, r)\right) + 1 + 2\sigma \tag{5.43}$$

$$\gamma_1(s) := \exp(\sigma\, r)\left(1 + r\, s\exp(\sigma\, r)\right) + 1 \tag{5.44}$$

Step $i = 2$: We define:

$$k_1(\xi_1) := -\left(\exp(\sigma\, r)\left(1 + r\,(1+\xi_1^2)\exp(\sigma\, r)\right) + 1 + 2\sigma\right)\xi_1 \tag{5.45}$$

$$\delta_1(\xi_1) := \left(\exp(\sigma\, r)\left(1 + r\,(1+3\xi_1^2)\exp(\sigma\, r)\right) + 1 + 2\sigma\right)\left(\exp(\sigma\, r)\left(1 + r\,(1+\xi_1^2)\exp(\sigma\, r)\right) + 2 + 2\sigma\right) \tag{5.46}$$

$$B_2(s) := \exp(\sigma\, r)\left(1 + r\,(1+s^2)\exp(\sigma\, r)\right) + 2 + 2\sigma \tag{5.47}$$

$$\begin{aligned}
\gamma_2(s) := {} & 2\exp(\sigma\, r)\left(\exp(\sigma\, r)\left(1 + r\,(1+4s^2\exp(2\sigma\, r))\exp(\sigma\, r)\right) + 2 + 2\sigma\right) \\
& + 4\exp(2\sigma\, r) rs\left(\exp(\sigma\, r)\left(1 + r\,(1+4s^2\exp(2\sigma\, r))\exp(\sigma\, r)\right) + 2 + 2\sigma\right)^2 + 1
\end{aligned} \tag{5.48}$$

$$\rho_1(s) := \exp(\sigma\, r)\left(1 + 2r\, s\exp(\sigma\, r)\right) + 1 \tag{5.49}$$

and

$$\mu_2(\xi_1, \xi_2) = \sigma + \gamma_2(p) + \gamma_1(p)\delta_1(\xi_1) + \frac{3}{4\sigma}\gamma_2^2(p) + \frac{3}{4\sigma}\gamma_1^2(p)\delta_1^2(\xi_1) + \frac{3}{4\sigma}\rho_1^2(p) \tag{5.50}$$

where

$$p := 1 + \xi_1^2 + \left|\xi_2 + \left(\exp(\sigma\, r)\left(1 + r\,(1+\xi_1^2)\exp(\sigma\, r)\right) + 1 + 2\sigma\right)\xi_1\right|^2 \tag{5.51}$$

The stabilizing feedback law is given by:

$$u(t) = -\mu_2(x_1(t), x_2(t))\left(x_2(t) + \left(\exp(\sigma\, r)\left(1 + r\,(1+x_1^2(t))\exp(\sigma\, r)\right) + 1 + 2\sigma\right)x_1(t)\right) \tag{5.52}$$

By virtue of Theorem 5.1, the functional

$$V(x) := \max_{\theta\in[-r,0]} \exp(2\sigma\,\theta)\left( x_1^2(\theta) + \left|x_2(\theta) + \left(\exp(\sigma\, r)\left(1 + r\,(1+x_1^2(\theta))\exp(\sigma\, r)\right) + 1 + 2\sigma\right)x_1(\theta)\right|^2 \right) \tag{5.53}$$

is a SRCLF which satisfies the small-control property and system (5.42) with (5.52) is URGAS. $\triangleleft$

## 6. Conclusions

In the present work we have showed how the well-known "Control Lyapunov Function (CLF)" methodology can be generalized to a broader class of nonlinear time-varying systems with both disturbance and control inputs, which include infinite-dimensional control systems described by retarded functional differential equations (RFDE). Necessary and sufficient conditions for the existence of stabilizing feedback are developed for the non-affine uncertain finite-dimensional case (1.1). Moreover, sufficient conditions, which guarantee that a given function is an Output Robust Control Lyapunov function for (1.1) are given. The case of uncertain control systems described by RFDEs of the form (1.3) is studied and special results are developed for the triangular case (1.4) of control systems described by RFDEs. It is shown that the construction of a stabilizing feedback law for (1.4) proceeds in parallel with the construction of a State Robust Control Lyapunov Functional. Moreover, sufficient conditions for the existence and design of a delay-free stabilizing feedback law are given. It is our belief that the present work can be used as a starting point for the discovery of necessary and sufficient Lyapunov-like conditions for the existence of stabilizing feedback for a wide class of infinite-dimensional control systems.

# Appendix

**Proof of Proposition 2.5:** Let $x:[t_0,t_{\max})\to\Re^n$ be the unique solution of (2.1) with initial condition $x(t_0)=x_0\in\Re^n$ corresponding to certain $d\in M_D$, where $t_{\max}\in(t_0,+\infty]$ is the maximal existence time of the solution. Notice that the mapping $[t_0,t_{\max})\ni t\to V(t,x(t))\in\Re^+$ is locally Lipschitz. By virtue of Lemma 2.4 and inequality (2.7) it follows that

$$\frac{d}{dt}\big(V(t,x(t))\big)\le-\rho\big(V(t,x(t))\big)+q(t)\text{, for all } t\in[t_0,t_{\max})\setminus N \tag{A1}$$

where $N\subset[t_0,t_{\max})$ is a set of zero Lebesgue measure.

<u>Case 1:</u> Hypothesis (Q2) holds.

Lemma 4.4 in [18] and inequality (A1) show that there exists $\sigma\in KL$ such that

$$V(t,x(t))\le\sigma\left(V(t_0,x_0),t-t_0\right)\text{, for all } t\in[t_0,t_{\max}) \tag{A2}$$

Inequality (A2) in conjunction with inequality (2.6) gives:

$$\mu(t)|x(t)|\le a_1^{-1}\big(\sigma\left(a_2(|x_0|),0\right)\big)\text{, for all } t\in[t_0,t_{\max}) \tag{A3}$$

Inequality (A3) and a standard contradiction argument show that (2.1) is RFC and that $t_{\max}=+\infty$. The fact that (2.1) is URGAOS follows from (A2) in conjunction with inequality (2.6). Particularly, we obtain:

$$|H(t,x(t))|\le a_1^{-1}\big(\sigma\left(a_2(|x_0|),t-t_0\right)\big)\text{, for all } t\ge t_0 \tag{A4}$$

Estimate (A4) implies that (2.1) is URGAOS.

<u>Case 2:</u> Hypotheses (Q1), (Q3) hold.

Lemma 3.2 in [10] and inequality (A1) show that there exists $\sigma\in KL$ such that

$$V(t,x(t))\le\sigma\left(V(t_0,x_0)+R,t-t_0\right)\text{, for all } t\in[t_0,t_{\max}) \tag{A5}$$

where $R:=\int_0^{+\infty} q(s)ds$ . Inequality (A5) in conjunction with inequality (2.6) gives:

$$\mu(t)|x(t)| \le a_1^{-1}\left(\sigma\left(a_2\left(\beta(t_0)|x_0|\right)+R,0\right)\right), \text{ for all } t\in[t_0,t_{\max}) \quad \text{(A6)}$$

Inequality (A6) and a standard contradiction argument show that (2.1) is RFC and that $t_{\max}=+\infty$ . The fact that (2.1) is RGAOS follows from (A6) in conjunction with inequality (2.6). Particularly, we obtain:

$$|H(t,x(t))| \le a_1^{-1}\left(\sigma\left(a_2\left(\beta(t_0)|x_0|\right)+R,t-t_0\right)\right), \text{ for all } t\ge t_0 \quad \text{(A7)}$$

Estimate (A7) shows that Property P2 (Uniform Output Attractivity on compact sets of initial data) holds of Definition 2.2 holds. Lemma 3.5 in [12] implies that (2.1) is RGAOS.

The proof is complete. ⊲

**Proof of Proposition 4.6:** Consider a solution $x(t)$ of (4.3) under hypotheses (Q1-5) corresponding to arbitrary $d\in M_D$ with initial condition $T_r(t_0)x=x_0\in C^1([-r,0];\Re^n)$. By virtue of Lemma 2.5 in [16], for every $T\in(t_0,t_{\max})$, the mapping $[t_0,T]\ni t\to V(t,T_r(t)x)$ is absolutely continuous. It follows from (4.8) and Lemma 4.4 that

$$\frac{d}{dt}\left(V(t,T_r(t)x)\right)\le -\rho\left(V(t,T_r(t)x)\right)+q(t) \text{ a.e. on } [t_0,t_{\max}) \quad \text{(A8)}$$

The previous differential inequality in conjunction with Lemma 3.2 in [10] shows that there exists $\sigma\in KL$ such that

$$V(t,T_r(t)x)\le\sigma\left(V(t_0,x_0)+R,t-t_0\right) \text{ for all } t\in[t_0,t_{\max}) \quad \text{(A9)}$$

where $R:=\int_0^{+\infty} q(s)ds$ . Using (4.7), (A9) and a standard contradiction argument we may show that $t_{\max}=+\infty$ . It follows from Lemma 2.6 in [16] that the solution $x(t)$ of (4.3) under hypotheses (Q1-5) corresponding to arbitrary $d\in M_D$ with arbitrary initial condition $T_r(t_0)x=x_0\in C^0([-r,0];\Re^n)$ satisfies (A9) for all $t\ge t_0$. Moreover, by virtue of (A9) and (4.7) we may establish that system (4.3) under hypotheses (Q1-5) is RFC. Notice that inequality (A9) in conjunction with (4.7) provide the estimate $\|H(t,T_r(t)x)\|_Y\le a_1^{-1}\left(\sigma\left(a_2(\beta(t_0)\|x_0\|_r)+R,t-t_0\right)\right)$ for all $t\ge t_0$, which establishes the fact that Property P2 (Uniform Output Attractivity on bounded sets of initial data) holds for system (4.3). It follows from Lemma 3.3 in [11] that system (4.3) is RGAOS. Furthermore, if $\beta(t)\equiv 1$ and $q(t)\equiv 0$, then inequality (A9) in conjunction with (4.7) provide the estimate $\|H(t,T_r(t)x)\|_Y\le a_1^{-1}\left(\sigma\left(a_2(\|x_0\|_r),t-t_0\right)\right)$ for all $t\ge t_0$, which establishes that system (4.3) is URGAOS.

The proof is complete. ⊲

**Proof of Lemma 5.3:** The fact that the functional $V$ as defined by (5.5) is Lipschitz on bounded sets of $C^0([-r,0];\Re^n)$ is a direct consequence of the fact that $Q\in C^1(\Re^n;\Re^+)$ (details are left to the reader). Consequently, as noticed in Section 4.II, we obtain the following simplification for the derivative $V^0(x;v)$ defined by (4.5) for all $(x,v)\in\Re^+\times C^0\left([-r,0];\Re^n\right)\times\Re^n$:

$$V^0(x;v)=\limsup_{h\to 0^+}\frac{V(E_h(x;v))-V(x)}{h}$$

Clearly, we have by virtue of (4.4) and (5.5):

$$V(E_h(x;v)) = \max\left\{ \max_{\theta\in[-r,-h]} \exp(2\sigma\theta)Q(x(\theta+h)),\ \max_{\theta\in[-h,0]} \exp(2\sigma\theta)Q(x(0)+(\theta+h)v) \right\}$$

$$= \max\left\{ \begin{array}{l} \exp(-2\sigma h) \max\limits_{s\in[-r+h,0]} \exp(2\sigma s)Q(x(s)), \\ \max\limits_{\theta\in[-h,0]} \exp(2\sigma\theta)\left[ Q(x(0)) + (\theta+h)\nabla Q(x(0))v + \int_0^{\theta+h} \left(\nabla Q(x(0)+sv) - \nabla Q(x(0))\right) v ds \right] \end{array} \right\} \quad \text{(A10)}$$

$$\leq \max\left\{ \exp(-2\sigma h)V(x),\ \max_{\theta\in[-h,0]} \exp(2\sigma\theta)\left[Q(x(0)) + (\theta+h)\nabla Q(x(0))v\right] + h|v| \max_{s\in[0,h]} \left|\nabla Q(x(0)+sv) - \nabla Q(x(0))\right| \right\}$$

If $Q(x(0)) < V(x)$ then there exists $h > 0$ such that $Q(x(0)) + s\nabla Q(x(0))v \leq \frac{1}{2}\left(V(x) + Q(x(0))\right)$ for all $s \in [0,h]$. Consequently, in this case we have from (A10) for $h > 0$ sufficiently small:

$$\begin{array}{l} h^{-1}(V(E_h(x;v)) - V(x)) \\ \leq \max\left\{ \dfrac{\exp(-2\sigma h) - 1}{h} V(x),\ \exp(-2\sigma h)\dfrac{1}{2h}\left(Q(x(0)) - V(x)\right) + |v| \max\limits_{s\in[0,h]} \left|\nabla Q(x(0)+sv) - \nabla Q(x(0))\right| \right\} \end{array}$$

The above inequality gives (5.6a) for the case $Q(x(0)) < V(x)$.

If $Q(x(0)) = V(x)$ and $\nabla Q(x(0))v > -2\sigma V(x)$ then it follows that $\max_{\theta\in[-h,0]} \exp(2\sigma\theta)\left[Q(x(0)) + (\theta+h)\nabla Q(x(0))v\right] = Q(x(0)) + h\nabla Q(x(0))v$ for $h > 0$ sufficiently small. Consequently, we obtain from (A10):

$$\begin{array}{l} h^{-1}(V(E_h(x;v)) - V(x)) \\ \leq \max\left\{ \dfrac{\exp(-2\sigma h) - 1}{h} V(x),\ \nabla Q(x(0))v + |v| \max\limits_{s\in[0,h]} \left|\nabla Q(x(0)+sv) - \nabla Q(x(0))\right| \right\} \end{array}$$

The above inequality gives (5.6b) for the case $Q(x(0)) = V(x)$ and $\nabla Q(x(0))v > -2\sigma V(x)$.

If $Q(x(0)) = V(x)$ and $\nabla Q(x(0))v \leq -2\sigma V(x)$ then it follows that $\max_{\theta\in[-h,0]} \exp(2\sigma\theta)\left[Q(x(0)) + (\theta+h)\nabla Q(x(0))v\right] = \exp(-2\sigma h)Q(x(0))$ for $h > 0$ sufficiently small. Consequently, we obtain from (A10):

$$h^{-1}(V(E_h(x;v)) - V(x)) \leq \frac{\exp(-2\sigma h) - 1}{h} V(x) + |v| \max_{s\in[0,h]} \left|\nabla Q(x(0)+sv) - \nabla Q(x(0))\right|$$

The above inequality gives (5.6b) for the case $Q(x(0)) = V(x)$ and $\nabla Q(x(0))v \leq -2\sigma V(x)$.

Thus inequality (5.6b) holds for the case $Q(x(0)) = V(x)$. The proof is complete. $\triangleleft$